\documentclass[10pt]{article}

\usepackage{mathrsfs}
\usepackage{amssymb}
\usepackage{amsmath}
\usepackage{amsfonts}
\usepackage{amstext}
\usepackage{amsthm}
\usepackage{graphicx}
\usepackage{subfig}
\usepackage{color}
\usepackage{indentfirst,latexsym,bm}
\usepackage{mathrsfs}
\usepackage{cite}
\usepackage[all]{xy}
\usepackage{newlfont}

\newcommand{\R}{\mathbb{R}}

\usepackage{enumerate}
\numberwithin{equation}{section}

\newtheorem{theorem}{Theorem}[section]

\newtheorem{lemma}[theorem]{Lemma}
\newtheorem{remark}[theorem]{Remark}

\newtheorem{proposition}[theorem]{Proposition}
\numberwithin{equation}{section}

\begin{document}

\title{Overdetermined elliptic problems in nontrivial exterior domains of the hyperbolic space
\thanks{The research of G. Dai is supported by the National Natural Science Foundation of China (Grants No.: 12371110). The research of P. Sicbaldi is supported by the FEDER-MINECO Grant PID2020-117868GB-I00 and the ``Maria de Maeztu'' Excellence Unit IMAG, reference CEX2020-001105-M, funded by MCINN/AEI/10.13039/501100011033/CEX2020-001105-M. The research of Y. Zhang is supported by the National Natural Science Foundation of China (Grants No.: 12301133).}}
\author{G. Dai
\thanks{Guowei Dai. School of Mathematical Sciences, Dalian University of Technology, Dalian, 116024, P.R. China. E-mail: daiguowei@dlut.edu.cn. } \,,
P. Sicbaldi
\thanks{Pieralberto Sicbaldi. Departamento de An\'alisis Matem\'atico \& Instituto de Matem\'aticas IMAG, Universidad de Granada, Avenida Fuentenueva, 18071 Granada, Spain \& Aix-Marseille Universit\'e - CNRS, Centrale Marseille - I2M, Marseille, France. E-mail: pieralberto@ugr.es. } \,,
Y. Zhang
\thanks{Yong Zhang. School of Mathematical Sciences, Jiangsu University, Zhenjiang, 212013, P.R. China. E-mail: zhangyong@ujs.edu.cn. }
}

\date{}
\maketitle

\renewcommand{\abstractname}{Abstract}

\begin{abstract}

We construct nontrivial unbounded domains $\Omega$ in the hyperbolic space $\mathbb{H}^N$, $N \in \{2,3,4\}$, bifurcating from the complement of a ball, such that the overdetermined elliptic problem
\begin{equation}
-\Delta_{\mathbb{H}^N} u+u-u^p=0\,\, \text{in}\,\,\Omega, \,\, u=0,\,\,\partial_\nu u=\text{const}\,\,\text{on}\,\,\partial\Omega\nonumber
\end{equation}
has a positive bounded solution in $C^{2,\alpha}\left(\Omega\right) \cap H^1\left(\Omega\right)$. We also give a condition under which this construction holds for larger dimensions $N$. This is linked to the Berestycki-Caffarelli-Nirenberg conjecture on overdetermined elliptic problems, and, as far as we know, is the first nontrivial example of solution to an overdetermined elliptic problem in the hyperbolic space. 
\end{abstract}

\\

\emph{AMS Subjection Classification(2020):} 35B32; 35N25; 51M10.


\section{Introduction and main results}

{\bf Introduction and state of the art.} Let us consider the following overdetermined elliptic problem:
\begin{equation}\label{ffunction}
\left\{
\begin{array}{ll}
\Delta u+f(u)=0\,\, &\text{in}\,\, \Omega,\\
u>0&\text{in}\,\, \Omega,\\
u=0 &\text{on}\,\, \partial \Omega,\\
\frac{\partial u}{\partial \nu}=\text{const} &\text{on}\,\, \partial \Omega,
\end{array}
\right.
\end{equation}
where $\Omega\subset \mathbb{R}^N$ is a given domain, $f$ is a Lipschitz function and $\nu$ is the unit outward normal vector about $\partial \Omega$. Serrin \cite{Serrin} proved that if $\Omega$ is bounded with boundary of class $C^2$ and $f$ is a $C^1$ function, the existence of a solution to (\ref{ffunction}) implies that $\Omega$ is a ball. The method used by Serrin to prove such result is universally known as the moving plane method and holds also when $f$ is only Lipschitz \cite{Pucci}. The case when $\Omega$ is supposed to be unbounded has been considered in 1997 by Berestycki, Caffarelli and Nirenberg \cite{BCN}. They obtained some rigidity results for epigraphs and stated the following conjecture:
\\ \\
\textbf{BCN Conjecture:} If $\Omega$ is a smooth domain such that $\mathbb{R}^N\setminus \overline{\Omega}$ is connected and problem (\ref{ffunction}) admits a bounded solution, then $\Omega$ is either a ball, a half-space, a generalized cylinder $B^k\times \mathbb{R}^{n-k}$ where $B^k$ is a ball in $\mathbb{R}^k$, or the complement of one of them.\\

The answer to the BNC conjecture, in its generality, is negative. In the class of domains diffeomorphic to a cylinder, for $N\geq3$ the second author \cite{Sicbaldi} found a periodic perturbation of the straight cylinder $B^{N-1}\times\mathbb{R}$ that supports a periodic solution to problem (\ref{ffunction}) with $f(u)=\lambda u, \lambda>0$. Such result holds also in dimension $N=2$ \cite{Schlenk} (but in this case is not a counterexample to the conjecture). Generalizations of this construction have been done in the Riemannian manifolds $\mathbb{S}^N\times \mathbb{R}$ and $\mathbb{H}^N\times \mathbb{R}$ \cite{DMS, Morabito}, and in the Euclidean case for general functions $f$ \cite{RSW22}. In the class of domains diffeomorphic to the half-space the BCN conjecture has been proved to be true in dimension 2 by Ros, Ruiz and the second author \cite{RRS} (under the assumption that the Neumann data is not zero), while new solutions in perturbations of the Bombieri-De Giorgi-Giusti epigraph have been built by Del Pino, Pacard and Wei in \cite{Del} in dimension $N\geq 9$ when $f(u)$ is of Allen-Cahn type (i.e. $f(u)=u-u^3$ or other nonlinearities with similar behavior). In the class of domains diffeomorphic to the exterior of a ball the BCN conjecture is true when $f$ is of Allen-Cahn type (see \cite{Reichel, RRS}), while new bounded solutions of \eqref{ffunction} in perturbations of the exterior of balls of large radius have been found by Ros, Ruiz and the second author \cite{RRS} for the Schr\"{o}dinger equation (i.e. $f(u)=u^p-u$). These new solutions in perturbations of the exterior of a ball produce automatically new solutions in perturbations of the complement of a straight cylinder just by adding one or more empty variables.

Two natural Riemannian manifolds where one can consider Problem (\ref{ffunction}) and hope to obtain rigidity results using the moving plane method and related tools are the round sphere $\mathbb{S}^N$ and the hyperbolic space $\mathbb{H}^N$ (and in fact these two manifolds are the only where one can hope this, as explained in the introduction of \cite{RSW2}). Hyperbolic geometry is very important in Physics: it is is closely related to the Minkowski space that arises from the Einstein's relativity and electromagnetic field theory \cite{Born0, Born1, Born2} and is also widely used in astrophysics, black hole theory and string theory.
Of course, when we consider a Riemannian manifold we replace the Laplacian operator in \eqref{ffunction} by the Laplace-Beltrami operator associated to the metric on the manifold.

Kumaresan and Prajapat \cite{Kumaresan} have proved that if $\Omega$ is a bounded domain of $\mathbb{H}^N$ or of an hemisphere of $\mathbb{S}^N$ and (\ref{ffunction}) admits a solution, then $\Omega$ must be a geodesic ball. Such result is the parallel of Serrin's theorem, and the proof is again based on the moving plane method. Nevertheless in the all sphere $\mathbb{S}^N$ this statement turns to be false even for domains diffeomorphic to a geodesic ball: in \cite{RSW2} Ruiz, Wu and the second author found a nontrivial contractible domain where problem (\ref{ffunction}) can be solved with $f(u)=u^3-u$ (up to a multiplicative constant). Such result shows that the nontrivial topology of $\mathbb{S}^N$ allows new phenomena and the situation in $\mathbb{S}^N$ can be very different from $\mathbb{R}^N$ concerning overdetermined elliptic problems.

In $\mathbb{H}^N$ the situation seems to be closer to the Euclidean one, but with some differences and the classes of possible shapes where one can hope to solve problems of the form \eqref{ffunction} are richer. For example, beside the classical domains (i.e. geodesic balls, half-spaces, the exterior of a geodesic ball, cylinders, the complement of a cylinder) there exists an other trivial shape where one can solve overdetermined problems: the horoball (see below its definition). After the pionner work \cite{Kumaresan}, overdetermined elliptic problems in $\mathbb{H}^N$ have been considered in a structural way in two very interesting papers: \cite{EM18} by Espinar ad Mao, and \cite{Espinar} by Espinar, Farina and Mazet. The hyperbolic space can be compactified by its ideal boundary $\partial_\infty \mathbb{H}^N$, and domains where a problem as \eqref{ffunction} admits solutions can be studied in terms of their trace on $\partial_\infty \mathbb{H}^N$. In order to give a rigorous definition of $\partial_\infty \mathbb{H}^N$, we say that two geodesics $\gamma_1(t)$ and $\gamma_2(t)$ are asymptotic if there exists a constant $c$ such that $d(\gamma_1(t),\gamma_2(t))<c$ for all $t\geq0$ where $d$ stands for the distance in $\mathbb{H}^N$ and $c$ is a positive constant. Asymptotic is an equivalence relation on the set of unit-speed half geodesics (i.e. geodesics considered only for $t\in[0,+\infty)$). Each equivalence class is called a point at infinity and we denote by $\partial_\infty \mathbb{H}^N$ the set of points at infinity. It is well known that $\partial_\infty \mathbb{H}^N$ is bijective to a unit sphere, i.e., $\partial_\infty \mathbb{H}^N \equiv \mathbb{S}^{N-1}$. Given a geodesic $\gamma$, the function
\[
B_\gamma(x) = \lim_{t \to +\infty} [d(\gamma(t),x) -t ]
\]
is called the Busemann function at $\gamma$. The level sets (respectively super-level sets) of such function are called horospheres (respectively horoballs). The ideal boundary of a horosphere (and a horoball) is just one point in $\partial_\infty \mathbb{H}^N$, given by $\gamma(+\infty)$.

The main result about $\mathbb{H}^N$ contained in \cite{EM18} is the following: if Problem \eqref{ffunction} has a bounded solution in a domain $\Omega \subset \mathbb{H}^N$ and the boundary at infinity of $\Omega$ stays in an equator of $\mathbb{S}^{N-1}=\partial_\infty \mathbb{H}^N$ then the domain must be symmetric with respect to the hyperplane whose boundary at infinity is such equator. In particular:
\begin{itemize}
\item If the boundary at infinity is done by exactly two different points (i.e. the domain has two ends) then the domain must be rotationally symmetric with respect to the axis given by the geodesic that joins the two points.
\item If the boundary at infinity is done by exactly one point, then the domain is a horoball.
\end{itemize}
The result for the case with two different points at infinity can be seen as the analogous to a result in $\mathbb{R}^N$ by Ros and the second author \cite{Ros}.

In \cite{Espinar} the authors consider Problem \eqref{ffunction} in the case of domains $\Omega$ of $\mathbb{H}^N$ diffeomorphic to the complement of a geodesic ball or the complement of a horoball. They assume the following hypothesis on the solution $u$ and on the function $f$:

(H) there exists a positive constant $C$ such that $u(p) \to C$ when $d(p,\partial \Omega) \to +\infty$, and $f$ is non-increasing in an interval $(C-\epsilon, C)$ for some (small) $\epsilon$.

Under such hypothesis, the two main results proved in \cite{Espinar} are:
\begin{itemize}
\item If $\Omega$ is the complement of a bounded open region with regular boundary, and Problem \eqref{ffunction} has a solution, then the domain must be the complement of a geodesic ball.
\item If $\Omega$ is the complement of an open domain with regular boundary and only one point at infinity, and Problem \eqref{ffunction} has a solution, then the domain must be the complement of a horoball.
\end{itemize}
In \cite{EM18, Espinar} the authors suggest implicitly the equivalent of the BCN Conjecture, that can be stated in the following problem: prove that, under some additional hypothesis on the function $f$ or on the geometry of the domain $\Omega$, the existence of a bounded solution to Problem \eqref{ffunction} implies that $\Omega$ must be a geodesic ball, a horoball, a cylinder (with base a geodesic ball or a horoball), a half-space, or the complement of one of them. In this sense, in \cite{EM18} Espinar and Mao prove the BCN Conjecture in $\mathbb{H}^N$ under the hypothesis that the boundary at infinity of the domain is empty or done by only one point, while in \cite{Espinar} the authors prove the BCN conjecture in $\mathbb{H}^2$ with the additional hypothesis (H).

The aim of this paper is to show that hypothesis (H) in \cite{Espinar} for the rigidity result on domains that are the complement of a geodesic ball cannot be eliminated. In fact, we construct nontrivial exterior domains which support a positive bounded solution of Problem (\ref{ffunction}) when $f$ is the nonlinear Schr\"{o}dinger function. That is,
here we consider the following problem
\begin{equation}\label{canshu}
\left\{
\begin{array}{ll}
-\Delta u+ u-u^p=0, \,\,u>0\,\, &\text{in}\,\, \Omega,\\
u=0 &\text{on}\,\, \partial \Omega,\\
\frac{\partial u}{\partial \nu}=\text{const} &\text{on}\,\, \partial \Omega,
\end{array}
\right.
\end{equation}
where $\Omega\subset\mathbb{H}^N$ with $N\geq 2$ and $1<p<(N+2)/(N-2)$ if $N>2$ and $1<p<+\infty$ if $N=2$. We construct nontrivial exterior domains $\Omega$, obtained as perturbations of the complement of a big geodesic ball, that support solutions to \eqref{canshu} converging to 0 when the distance to $\partial \Omega$ tends to infinity.

Let us mention that the study of the nonlinear Schr\"{o}dinger equation in the hyperbolic space
is motivated by the Grushin operators \cite{Beckner}, Hardy-Sobolev-Mazy's equations \cite{Castorina} and quantum mechanics\cite{Banica}. There is a large amount of literature dealing with the study of existence, uniqueness and
qualitative properties of solutions to the nonlinear Schr\"{o}dinger equation in the hyperbolic space, see for example \cite{Banica, Mancini, Morabito0, Morabito} and references therein.\\ 

{\bf Statement of the results.} In order to state precisely our result we need some notation. Fixing the origin $O$ in $\mathbb{H}^N$, we use the exponential map of $\mathbb{H}^N$ centered at $O$, $\mbox{exp}_O: \mathbb{R}^N \to \mathbb{H}^N$ and given any continuous function $v: \mathbb{S}^{N-1} \to \left(0, \infty \right)$, we define the domain
\[  B_v = \mbox{exp}_O \left (  \left\{x\in\mathbb{H}^{N}:0\leq|x|<v\left(\frac{x}{|x|}\right)\right\} \right )\,.
\]
and let $B_v^c$ be its complement in $\mathbb{H}^N$, where $|x|$ is the geodesic distance of $x$ to the origin (we will write $B_R$ for the ball of radius $R$, i.e. $B_v$ with $v(x)\equiv R$). We shall also use the coordinates in $\mathbb{H}^N$ given by the exponential map composed with polar coordinates in $\R^N$. In other words, we write:
\begin{equation} \label{coord}  \begin{array}{c} X: \left[0, \infty \right) \times \mathbb{S}^{N-1} \to \mathbb{H}^N \\
X(r, \theta) = \mbox{exp}_S(r\,  \theta). \end{array} \end{equation}
Let $G$ be a group of isometries acting on $\mathbb{S}^{N-1}$. We say that
$\Omega\subset \mathbb{H}^N$ is $G$-symmetric if, working in coordinates:
\[
(r, \theta) \in \Omega \Rightarrow (r, g(\theta)) \in \Omega
\]
for any $g \in G$. Through the paper we take $\alpha \in (0, 1)$ fixed. Define
\begin{equation}
H_{0,G}^1(\Omega)=\left\{u\in H_0^1(\Omega):u(r,\theta)=u (r, g(\theta)),\, \forall g\in G\right\},\nonumber
\end{equation}
\begin{equation}
C_G^{k,\alpha}(\Omega)=\left\{u\in C^{k,\alpha}(\Omega):u(r,\theta)=u (r, g(\theta)) ,\, \forall g\in G\right\} \nonumber
\end{equation}
\begin{equation}
C_G^{k,\alpha}(\mathbb{S}^{N-1})=\left\{u\in C^{k,\alpha}(\mathbb{S}^{N-1}):u=u \circ g, \, \forall g\in G\right\} \, . \nonumber
\end{equation}
Let $C_{G,m}^{k,\alpha}\left(\mathbb{S}^{N-1}\right)$ be the set of functions in $C_G^{k,\alpha}\left(\mathbb{S}^{N-1}\right)$ whose mean
is $0$. We assume that the group $G$ satisfies the following fundamental hypothesis:\\

(G1) Denoting by $\left\{\mu_{i_k}
\right\}_{k\in \mathbb{N}}$ the eigenvalues of $-\Delta_{\mathbb{S}^{N-1}}$
restricted to $G$-symmetric functions and by $m_k$ their multiplicities, we require $m_1$ odd and $i_1 >\left(2-N+\sqrt{(N-2)^2+\frac{16}{9}(N+2)(N-1)}\right)/2$.\\

We remark that when $N=2, 3$ or $4$ the condition is equivalent to $i_1 \geq 2$. So, in dimension $2$ any dihedral group satisfies (G1); while, in dimension $3$, $G$ can be taken as the group of isometries of the tetrahedron, the octahedron or the icosahedron (see \cite[Remark 2.2]{RRS}). In dimension $4$, $G$ can be the symmetry group of rotations of the hyper-icosahedron, for which $i_1=12$ and $m_1=1$ (see \cite[Appendix 7.3]{Ruiz}). For higher dimensions, we do not know if there exist groups satisfying (G1).

By the Krasnosel'skii Local Bifurcation Theorem \cite{Krasnoselskii}, we obtain the following result.
\begin{theorem}\label{teo11}
Let $N\in \mathbb{N}$ with $N\geq 2$, $1 < p < (N+2)/(N-2)$
($p > 1$ if $N = 2$). Let $G$ be a group
of symmetries of $\mathbb{S}^{N-1}$ satisfying (G1). Then there exist $R_* > 0$, a sequence
of nonzero functions $v_n\in C_{G,m}^{2,\alpha}\left(\mathbb{S}^{N-1}\right)$ converging to $0$, and a sequence of positive real
numbers $R_n$ converging to $R_*$ such that the overdetermined problem
\begin{equation}
\left\{
\begin{array}{ll}
-\Delta u+ u-u^p=0, \,\,u>0\,\, &\text{in}\,\, B_{R_n\left(1+v_n\right)}^c,\\
u=0 &\text{on}\,\, \partial B_{R_n\left(1+v_n\right)},\\
\frac{\partial u}{\partial \nu}=\text{const} &\text{on}\,\, \partial B_{R_n\left(1+v_n\right)}
\end{array}
\right.\nonumber
\end{equation}
has a positive bounded solution $u\in C_G^{2,\alpha}\left(B_{R_n\left(1+v_n\right)}^c\right)\cap H_{0,G}^1\left(B_{R_n\left(1+v_n\right)}^c\right)$.
\end{theorem}

Theorem \ref{teo11} is exactly the counterpart in $\mathbb{H}^N$ of the homologous result in $\mathbb{R}^N$ \cite[Theorem 1.1]{RRS}, but due to the hyperbolic space, we need a stronger condition for $i_1$ than in the Euclidean space where $i_1\geq2$ \cite{RRS}. 
We will follow the same strategy of \cite{RRS}, but we will have to treat the aspects related to the  hyperbolic metric.\\

{\bf Organization of the paper.} The outline of the paper is as follows. In Section 2 we study the existence, uniqueness, non-degeneracy and asymptotic behavior of the radial solutions of the Dirichlet problem with nonlinear Schr\"{o}dinger
equation in the exterior domain $B_R^c\subset\mathbb{H}^N$ for $N\geq 2$. 
After a change of scale and considering the parameter $\lambda = 1/R^2$, in Section 3 we introduce and study two quadratic forms $Q_{\lambda}(\psi)$ and $\widetilde{Q}_{\lambda}(\psi)$. In section 4, after solving the Dirichlet problem in domains of the form $B_{R(1+v)}^c$, we define the main operator of our problem, that turns to be the operator associating to $v$ the normal derivative at the boundary of the solution to the Dirichlet problem in $B_{R(1+v)}^c$, and we study the eigenvalues $\sigma_k\left(H_\lambda\right)$
of its linearization $H_{\lambda}$. 
This will allow us to use the classical bifurcation theory and in Section 5 we provide the proof of our main result.

\section{The Dirichlet problem in the exterior of a ball in $\mathbb{H}^N$}

For convenience of the reader, let's first recall the fundamental definition of the hyperbolic space.
For $x,y\in \mathbb{R}^{N+1}$ with $x=\left(x_0,x_1,\ldots,x_N\right)$ and $y=\left(y_0,y_1,\ldots,y_N\right)$, define the Lorentz inner product on $\mathbb{R}^{N+1}$ by
\begin{equation}
\left[x,y\right]=x_0y_0-\sum_{i=1}^Nx_iy_i.\nonumber
\end{equation}
The branch of
hyperboloid given by
\begin{equation}
\left\{x\in \mathbb{R}^{N+1}:[x,x]=1,x_0>0\right\}\nonumber
\end{equation}
equipped with the metric induced by the Lorentz inner product on $\mathbb{R}^{N+1}$ is called the hyperbolic space and denoted by $\mathbb{H}^N$.
Fix the point $O=(1,0,...,0)$ as the origin of $\mathbb{H}^N$. We will use spherical coordinates $(r,\omega)\in [0,+\infty)\times \mathbb{S}^{N-1}$ centered at $O$. This means that if $(x_0,x') \in \mathbb{H}^N$, $x'= (x_1, x_2, ..., x_N)$, we denote by $r$ the distance of $(x_0,x')$ to $O$ and we will write
\[
(x_0,x') = \left(\cosh r,\omega\sinh r\right)\,.
\]
We see that
\begin{equation}
\text{d}x_0=\sinh r\,\text{d}r,\,\,\text{d}x'=\omega\cosh r \,\text{d}r+\sinh r\,\text{d}\omega.\nonumber
\end{equation}
Then the metric induced on $\mathbb{H}^N$ by the Lorenzian one on $\mathbb{R}^{N+1}$ is
\begin{equation}
\text{d}l^2=-\text{d}x_0^2+\text{d}(x')^2=\text{d}r^2+\sinh^2r\,\text{d}\omega^2.\nonumber
\end{equation}
The Laplace-Beltrami operator $\Delta$ on $\mathbb{H}^N$ can then be written as
\begin{equation}
\partial_r^2+(N-1)\frac{\cosh r}{\sinh r}\partial_r+\frac{1}{\sinh^2 r}\Delta_{\mathbb{S}^{N-1}}\nonumber
\end{equation}
where $\Delta_{\mathbb{S}^{N-1}}$ is the Laplace-Beltrami operator on $\mathbb{S}^{N-1}$.
We refer the reader to the monographs \cite{Helgason, Terras} for more details on the hyperbolic space.\\

{\bf The radial solution to the Dirichlet problem in the exterior of a ball.} Consider the following problem
\begin{equation}\label{equationBR}
\left\{
\begin{array}{ll}
-\Delta w+w-w^p=0\,\, &\text{in}\,\, B_R^c,\\
w=0 &\text{on}\,\, \partial B_R,
\end{array}
\right.
\end{equation}
where $p>1$. Our first task will be to show that there exists a unique radial solution $w$ of \eqref{equationBR}, belonging to $H^1(B_R^c)$, non degenerate among radial functions and with Morse index equal to 1. For this we will use \cite{Morabito0}, where Morabito studies the existence, uniqueness and nondegeneracy of radial solutions of a more general equation in annular domains of more general Riemannian manifolds. We observe that if $w$ is radial and belongs to $H^1(B_R^c)$ then it satisfies the following ODE:
\begin{equation}\label{eigenvalueonbal001}
\left\{
\begin{array}{ll}
w''+(N-1)\frac{C(r)}{S(r)}w'+w^p-w=0,\,\,\, r\in(R,+\infty),\\
w(R)=0=\lim_{r\rightarrow+\infty}w(r),
\end{array}
\right.
\end{equation}
where
\begin{equation}
S(r)=\sinh (r)\quad \text{and}\quad C(r)=\cosh (r).\nonumber
\end{equation}
In the notation of \cite{Morabito0} this corresponds to $R_1=R$, $R_2 = +\infty$, $\nu = N-1$ and $V(r) =-1$. According to \cite{Morabito0}, in order to obtain the existence, the uniqueness and the nondegeneracy of a solution $w$ of \eqref{equationBR} we just need to study the derivative of the following function:
\begin{equation}
G(r)=-S^\beta(r)+\alpha S^{\beta-2}(r)\left[(\alpha+2-N)S'^2(r)-S''(r)S(r)\right],\nonumber
\end{equation}
where
\begin{equation}
\alpha=\frac{2(N-1)}{p+3},\,\,\beta=\alpha(p-1).\nonumber
\end{equation}
A straightforward computation (that we do in detail) gives us the following:\\
\begin{lemma}\label{lemma21}
Assume that $1 < p <(N+2)/(N-2)$
($p > 1$ if $N = 2$). Then either:
\begin{itemize}
\item $G'(r)< 0$ for $r \in (R,+\infty)$, or
\item $G'(R) > 0$ and $G'(r)$ changes sign only once on $(R,+\infty)$.
\end{itemize}
In both cases $\lim \inf_{r\rightarrow+\infty} G'(r) < 0$.
\end{lemma}
\textbf{Proof.} By the definition of $G$, we have that
\begin{eqnarray}
G'(r)&=&-\beta S^{\beta-1}(r)S'(r)+\alpha(\beta-2) S^{\beta-3}(r)S'(r)\left[(\alpha+2-N)S'^2(r)-S^2(r)\right]\nonumber\\
& &+2\alpha S^{\beta-2}(r)\left[(\alpha+2-N)S'(r)S(r)-S'(r)S(r)\right]\nonumber\\
&=&S^{\beta-1}(r)S'(r)\left[2\alpha(\alpha+1-N)-\beta-\alpha(\beta-2)+\alpha(\beta-2)(\alpha+2-N) S^{-2}(r)S'^2(r)\right].\nonumber
\end{eqnarray}
Since $1 < p < (N+2)/(N-2)$
($p > 1$ if $N = 2$), we derive that
\begin{eqnarray}
\beta-2=2\frac{(N-2)p-(N+2)}{p+3}<0.\nonumber
\end{eqnarray}
When $N=2$, we get that
\begin{eqnarray}
G'(r)&<&S^{\beta-1}(r)S'(r)\left[2\alpha(\alpha+1-N)-\beta-\alpha(\beta-2)+\alpha(\beta-2)(\alpha+2-N)\right]\nonumber\\
&=&S^{\beta-1}(r)S'(r)\beta(\alpha(\alpha-1)-1)=-S^{\beta-1}(r)S'(r)\frac{2(p-1)}{p+3}\left(1+\frac{2(p+1)}{(p+3)^2}\right)\nonumber\\
&<&0.\nonumber
\end{eqnarray}
Next we consider the case of $N\geq3$.

When $N\geq3$, we see that
\begin{eqnarray}
\alpha+2-N=\frac{(2-N)(p+1)+2}{p+3}\leq\frac{2-(p+1)}{p+3}<0.\nonumber
\end{eqnarray}
Define
\begin{eqnarray}
F(r)&=&2\alpha(\alpha+1-N)-\beta-\alpha(\beta-2)+\alpha(\beta-2)(\alpha+2-N) S^{-2}(r)S'^2(r).\nonumber
\end{eqnarray}
Then we have that
\begin{eqnarray}
F'(r)&=&2\alpha(\beta-2)(\alpha+2-N) \frac{C(r)}{S(r)}\left(1-\frac{C^2(r)}{S^2(r)}\right)<0.\nonumber
\end{eqnarray}
If $G'(R)\leq0$, we see that $F(R)\leq0$. Since $F$ is strictly decreasing, $F(r)<0$ on $(R,+\infty)$.
Thus, we have that $G'(r)<0$ on $(R,+\infty)$. If $G'(R)>0$, we have that $F(R)>0$. We note that
\begin{eqnarray}
\lim_{r\rightarrow+\infty}F(r)&=&\beta(\alpha(\alpha+1-N)-1)=-\frac{2(N-1)(p-1)}{p+3}\left(\frac{2(N-1)^2(p+1)}{(p+3)^2}+1\right)<0.\nonumber
\end{eqnarray}
So $F$ has at least one zero on $(R,+\infty)$. Since $F$ is strictly decreasing, $F$ has exactly one zero on $(R,+\infty)$.
Therefore, $G'(r)$ changes sign only once on $(R,+\infty)$.

Finally, we study $\lim \inf_{r\rightarrow+\infty} G'(r)$.
We first study the $\lim_{r\rightarrow+\infty} F(r)$.
When $N=2$, we have that
\begin{eqnarray}
\lim_{r\rightarrow+\infty} F(r)=-\frac{2(p-1)}{p+3}\left(1+\frac{2(p+1)}{(p+3)^2}\right)<0.\nonumber
\end{eqnarray}
And for $N\geq3$, we have that
\begin{eqnarray}
\lim_{r\rightarrow+\infty} F(r)=-\frac{2(N-1)(p-1)}{p+3}\left(\frac{2(N-1)^2(p+1)}{(p+3)^2}+1\right)<0.\nonumber
\end{eqnarray}
Note that
\begin{eqnarray}
G'(r)=S^{\beta}(r)\frac{S'(r)}{S(r)}F(r)\,\,\text{and}\,\,\lim_{r\rightarrow+\infty}S^{\beta}(r)\frac{S'(r)}{S(r)}=+\infty.\nonumber
\end{eqnarray}
So, we reach that
$\lim \inf_{r\rightarrow+\infty} G'(r)=-\infty<0$, which is the desired conclusion.\qed
\\

In view of Lemma \ref{lemma21}, by \cite[Theorem 2.1]{Morabito0}, the problem (\ref{eigenvalueonbal001}) has a unique positive solution $w_R$. Considering $w_R$ as radial solution of the problem in $B_R^c$, we have that it is non-degenerate in the space of $H_{0,r}^1\left(B_R^c\right)$, where
\begin{equation}
H_{0,r}^1\left(B_R^c\right)=\left\{z\in H_0^1\left(B_R^c\right):z(x)=z(\vert x\vert)\,\,\text{a.e.}\,\,x\in B_R^c\right\}.\nonumber
\end{equation}
We next derive a uniform decay rate of $w_R$.
\begin{lemma} \label{lemma22} Let $\gamma=\left(N-1+\sqrt{4+(N-1)^2}\right)/2$. For any $\varepsilon\in(0,\gamma)$, then there exists $M>0$ such that $w_R(r)\leq Me^{-r(\gamma-\varepsilon)}$ for $r$ large enough and $\lim_{r\rightarrow+\infty}w_R'/w_R=-\gamma$.
\end{lemma}
\noindent\textbf{Proof.}
Integrating the equation in (\ref{eigenvalueonbal001}) from $r$ to $+\infty$ and using the fact of $w_R\in H^1(B_R^c)$, we obtain that
\begin{eqnarray}
S^{N-1}(r)w_R'(r)=\int_r^{+\infty}S^{N-1}(\tau)\left(w_R^p(\tau)-w_R(\tau)\right)\,\text{d}\tau.\nonumber
\end{eqnarray}
Since $w_R(r)>0$ for any $r>R$, $\lim_{r\rightarrow+\infty}w_R(r)=0$ and $p>1$, we can take $r$ large enough such that
\begin{eqnarray}
\int_r^{+\infty}S^{N-1}(\tau)\left(w_R^p(\tau)-w_R(\tau)\right)\,\text{d}\tau\leq0.\nonumber
\end{eqnarray}
It follows that
$w_R'(r)\leq0$ for $r$ large enough.
Further, we have that
\begin{eqnarray}
w_R''(r)=-(N-1)\frac{S'(r)}{S(r)}w_R'(r)+w_R(r)-w_R^p(r)>0\nonumber
\end{eqnarray}
for $r$ large enough.
It follows that
\begin{eqnarray}
w_R'(r)<0\nonumber
\end{eqnarray}
for $r$ large enough.
Indeed, if there exists some $r_0$ large enough such that $w_R'\left(r_0\right)=0$, then we find that
\begin{eqnarray}
0\geq\lim_{\varepsilon\rightarrow0^+}\frac{w_R'\left(r_0+\varepsilon\right)}{\varepsilon}=w_R''\left(r_0\right)>0,\nonumber
\end{eqnarray}
which is a contradiction.
We take $R_0>R$ such that $w_R'(r)<0$ for $r>R_0$.

For any $r>R_0$, set $z=-w_R'/w_R$.
Observe that $z>0$. Let us prove that $z$ is also bounded from above.

We have:
\begin{equation}
z'=-\frac{w_R''}{w_R}+\frac{w_R'^2}{w_R^2}=z^2-(N-1)\frac{S'}{S}z+w_R^{p-1}-1.\nonumber
\end{equation}
Let us write:
\begin{equation}
z^2-(N-1)\frac{S'}{S}z+w_R^{p-1}-1=\frac{z^2}{2}+h(z),\nonumber
\end{equation}
where
\begin{equation}
h(z):=\frac{z^2}{2}-(N-1)\frac{S'}{S}z+w_R^{p-1}-1.\nonumber
\end{equation}
We observe that
\begin{equation}
\lim_{r\rightarrow+\infty}\left((N-1)^2\frac{S'^2}{S^2}+2\left(1-w_R^{p-1}\right)\right)=2+(N-1)^2>0.\nonumber
\end{equation}
So there exists $R_1\geq R_0$ such that the discriminant of the second order polynome $h$ is positive if $r\geq R_1$.
It follows that, if $r>R_1$, the polynome $h$ has exactly two zeros, and the larger one is
\begin{equation}
z_0:= (N-1)\frac{S'}{S}+\sqrt{(N-1)^2\frac{S'^2}{S^2}+2\left(1-w_R^{p-1}\right)}\,,\nonumber
\end{equation}
which is clearly positive. Suppose now that $z\geq z_0$ for $r$ in a certain interval of $(R_0,+\infty)$. In this case we have $h(z)\geq0$.
It follows that
\begin{equation}
z'\geq \frac{z^2}{2}\,.\nonumber
\end{equation}
So $z$ is stricly increasing and
\begin{equation}\label{integratingz}
\left(\frac{1}{z(r)}\right)'\leq -\frac{1}{2}.
\end{equation}
Note that
\begin{equation}
\lim_{r\rightarrow+\infty} z_0=N-1+\sqrt{2+(N-1)^2}\,.\nonumber
\end{equation}
Thus, there exists $R_2\geq R_1$ such that
\begin{equation}
z_0<N+\sqrt{2+(N-1)^2}= :N^*\nonumber
\end{equation}
for $r\geq R_2$.

So, now suppose by contradiction that $z$ is bigger than $N^*$ at some $r \geq R_2$. Since we can choose $R_2$ larger, let us suppose that $z(R_2)\geq N^*$. Then, since $z$ is strictly increasing in intervals where $z \geq z_0$, 
we have that $z$ is strictly increasing for $r \geq R_2$, so $z \geq N^*$ for all $r \geq R_2$. 
%
Integrating (\ref{integratingz}) from $R_2$ to $r$, we find that
\begin{equation}
\int_{R_2}^{r}\left(\frac{1}{z(t)}\right)'\,\text{d}t\leq \frac{R_2-r}{2},\nonumber
\end{equation}
which means
\begin{equation}
\frac{1}{z(r)}\leq \frac{2-\left(r-R_2\right)z\left(R_2\right)}{2z\left(R_2\right)}. \nonumber
\end{equation}
Since $z$ is positive for $r\geq R_2$, then it's obvious that
\begin{equation}
z(r)\geq \frac{2z\left(R_2\right)}{2-\left(r-R_2\right)z\left(R_2\right)}, \nonumber
\end{equation}
which implies that $z(r)$ would blow up as $r\rightarrow \left(2+R_2z\left(R_2\right)\right)/z\left(R_2\right)$, that is a value in $(R_2, +\infty)$. Due to the fact that $z(r)$ is differentiable for any $r>R_0$, this behavior is impossible and we have a contradiction with our assumption $z(R_2) \geq N^*$.

It follows at once that $0 < z < N^*$ for $r \geq R_2$ and then
\begin{equation}
0\leq\liminf_{r\rightarrow+\infty}z\leq\limsup_{r\rightarrow+\infty}z\leq N^*\,. \nonumber
\end{equation}
Now, using the argument of \cite[Lemma 3.4]{Mancini} with obvious changes, we can show that
\begin{equation}
\liminf_{r\rightarrow+\infty}z=\limsup_{r\rightarrow+\infty}z\, , \nonumber
\end{equation}
and so $z$ has a limit for $r \to +\infty$.
Let ${\gamma}=\lim_{r\rightarrow+\infty}z$.
Then by L'H\^{o}spital's rule, we have that
\begin{equation}
{\gamma}^2=\lim_{r\rightarrow+\infty}z^2=\lim_{r\rightarrow+\infty}\frac{w_R''}{w_R}
=\lim_{r\rightarrow+\infty}\left((N-1)\frac{S'}{S}z+1-w_R^{p-1}\right)=1+(N-1){\gamma}.\nonumber
\end{equation}
It follows that
\begin{equation}
{\gamma}=\frac{N-1+\sqrt{4+(N-1)^2}}{2}>0.\nonumber
\end{equation}
For any $\varepsilon\in(0,\gamma)$ and $r$ large enough, we have that
\begin{equation}
-\frac{w_R'}{w_R}\geq\gamma-\varepsilon.\nonumber
\end{equation}
It follows that
\begin{equation}
\ln\frac{1}{w_R}\geq(\gamma-\varepsilon)r+C\nonumber
\end{equation}
for some constant $C$ and $r$ large enough.
Therefore, for some positive constant $M$, we obtain that
\begin{equation}
w_R\leq Me^{-(\gamma-\varepsilon)r},\nonumber
\end{equation}
which is the desired conclusion.\qed
\\

{\bf Remark.} In the previous proof we combine the methods of proving decay rate in $\mathbb{R}^N$ of \cite{Peletier} and $\mathbb{H}^N$ \cite{Mancini}.
In Euclidean space, if $w_R$ is bounded, $w_R/r\rightarrow0$ which implies the existence of $\lim_{r\rightarrow+\infty}z$.
While, in the hyperbolic space, since $S'/S\rightarrow1$, the L'H\^{o}spital rule used in the Euclidean space can not be directly used here.
This is also an essential difference between the hyperbolic space and the Euclidean space.
Of course, the existence of $\lim_{r\rightarrow+\infty}z$ can also be obtained by using the method developed in \cite{Mancini}.
Here we combine two methods to display the differences between two spaces.\\

Furthermore, the radial solution $w_R$ also has the following property.\\
\begin{lemma} \label{lemma23} The function $w_R$ increases in the radius up to a certain maximum, and then
it decreases and converges to $0$ at infinity.
\end{lemma}
\textbf{Proof.} Define
\begin{equation}
E(r)=\frac{w_R'^2(r)}{2}+\int_0^{w_R(r)}\left(s^p-s\right)\,\text{d}s,\nonumber
\end{equation}
then
\begin{equation}
E'(r)=-(N-1)w_R'^2(r)\frac{S'(r)}{S(r)}\leq0.\nonumber
\end{equation}
Since $w_R(R)=0=\lim_{r\rightarrow+\infty}w_R(r)$, $w_R$ has at least one critical point.
Let $c$ be any critical point of $w_R$.
One see that
\begin{equation}
E(c)=\int_0^{w_R(c)}\left(s^p-s\right)\,\text{d}s=\frac{w_R^{p+1}(c)}{p+1}-\frac{w_R^2(c)}{2}.\nonumber
\end{equation}
In view of Lemma \ref{lemma22}, one may see that $\lim_{r\rightarrow+\infty}E(r)=0$.
It follows that $E(r)\geq0$. In particular, $E(c)\geq0$.
It follows that
\begin{equation}
w_R(c)\geq \left(\frac{p+1}{2}\right)^{\frac{1}{p-1}}>1.\nonumber
\end{equation}
So, the critical values of $w_R$ is greater than $1$. We claim that the maximum point of $w_R$ is unique as in the following Figure 1.
\begin{figure}[ht]
\centering
\includegraphics[width=1.0\textwidth]{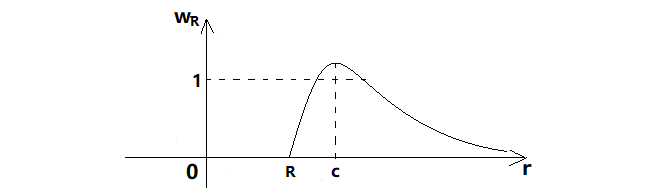}
\caption{The behavior of the radial solution $w_R$.}
\end{figure}

Suppose, by contradiction, that
there exist $c$, $d$ with $c<d$ such that $w_R(c)=w_R(d)=\max_{[R,+\infty)} w_R(r)$.
Since $w_R'(c)=w_R'(d)=0$, we have that
\begin{equation}
E(c)=E(d).\nonumber
\end{equation}
The monotonicity of $E$ implies that $E(r)\equiv E(c)$ on $[c,d]$.
So $E'(r)\equiv0$ on $[c,d]$. Further, $w_R'(r)\equiv0$ on $[c,d]$.
It follows from (\ref{eigenvalueonbal001}) that $w_R(r)\equiv1$ on $[c,d]$, which contradicts the fact of any critical values of $w_R$ is greater than $1$.
Finally, we show that $w_R$ has exactly one critical point.
If it was not, there would exist a point $e$ such that $w_R'(e)=0$, $w_R''(e)\geq0$.
By \eqref{eigenvalueonbal001}, we have that
\begin{equation}
\left(w_R^p-w_R\right)\left(e\right)\leq0.\nonumber
\end{equation}
This implies that $w_R(e)\leq 1$, which again contradicts the fact of any critical values of $w_R$ is greater than $1$ (see the following Figure 2).
\begin{figure}[ht]
\centering
\includegraphics[width=1.0\textwidth]{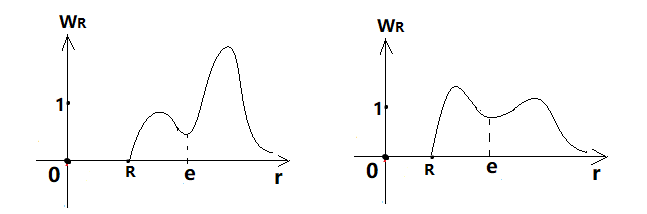}
\caption{The cases for multiple critical points of $w_R$. }
\end{figure}

Therefore, $w_R$ increases in the radius up to a certain maximum, and then
it decreases and converges to $0$ at infinity.
\qed\\

In addition, we get the exact Morse index of $w_R$ as follows.
\begin{lemma}
The function $w_R$ has Morse index equal to $1$.
\end{lemma}
\textbf{Proof.} For $u\in H_{0,r}^1\left(B_R^c\right)$, define the functional
\begin{equation}
I(u)=\frac{1}{2}\int_{B_R^c}\left(\vert \nabla u\vert^2+u^2\right)\,\text{d}x-\frac{1}{p+1}\int_{B_R^c} u_+^{p+1}\,\text{d}x,\nonumber
\end{equation}
where $u_+=\max\{u,0\}$.
The unique solution $w_R$ is a critical point of $I$.
It is standard to verify that $I$ has mountain pass geometry and satisfies the Palais-Smale compactness condition.
Since $w_R$ is the unique non-zero critical point of $I(u)$, $w_R$ is the mountain pass critical
point. So the Morse index $i\left(I,w_R\right) \leq 1$ \cite{Hofer} where $i\left(I,w_R\right)$ is defined by
\begin{eqnarray}
\max\left\{\text{dim} H:H\subseteq H_{0,r}^1\left(B_R^c\right), I''\left(w_R\right)(h,h)<0,\forall h\in H\setminus\{0\}\right\}.\nonumber
\end{eqnarray}
Observe that
\begin{eqnarray}
I''\left(w_R\right)(h,h)=\int_{B_R^c}\left(\vert \nabla h\vert^2+h^2-pw_R^{p-1}h^2\right)\,\text{d}x.\nonumber
\end{eqnarray}
Choosing $h = w_R$ and using $I'\left(w_R\right)w_R= 0$, we get
\begin{eqnarray}
I''\left(w_R\right)\left(w_R,\right)=-(p-1)\int_{B_R^c} w_R^{p+1}\,\text{d}x<0.\nonumber
\end{eqnarray}
So we get $i\left(I,w_R\right) = 1$.
\qed\\

For any $x\in \overline{B}_1^c\subset\mathbb{H}^N$, let $u(x)=w(Rx)$. Then $u$ satisfies
\begin{equation}\label{eigenvalueonbal01}
\left\{
\begin{array}{ll}
-\lambda\Delta u+u-u^p=0\,\, &\text{in}\,\, B_1^c,\\
u=0 &\text{on}\,\, \partial B_1,
\end{array}
\right.
\end{equation}
where $\lambda=1/R^2$.
From the properties of $w_R$, we have the following conclusion.
\begin{proposition}\label{prop21}
For any $\lambda > 0$, problem (\ref{eigenvalueonbal01}) possesses a unique radially symmetric $C^\infty$ solution $u_\lambda$, which increases in the radius up to a certain maximum, and then
it decreases and converges to $0$ at infinity. Moreover, $u_\lambda$ is non-degenerate in the space $H_{0,r}^1\left(B_1^c\right)$ and has Morse index $1$.
\end{proposition}
We use $\dot{u}_\lambda$ to denote its derivative with respect to $\lambda$. Then from the relations of $w$ and $u_\lambda$ we have that
\begin{equation}\label{uderivative}
\dot{u}_\lambda(r)=-u_\lambda'(r)\frac{r}{2R}\lambda^{-\frac{3}{2}}=-u_\lambda'(r)\frac{r}{2\lambda},
\end{equation}
which will be used later.
We regard $\left(u_\lambda,B_1^c\right)$ as the trivial solution pair of problem (\ref{canshu}).\\

\section{A preliminar tool: study of two quadratic forms}
The main goal of the present paper is to find a nontrivial exterior domain which supports a positive bounded solution of problem (\ref{canshu}) by deforming the complement of the ball.
To this end, we divide this process into two steps. The first step aims to find a solution of the problem (\ref{eigenvalueonbal01}) on a perturbation domain. After that we impose the Neumann condition on this solution by solving an operator equation.

In order to solve (\ref{eigenvalueonbal01}) on a perturbation domain, we consider its linearized eigenvalue problem near the radial solution $u_\lambda$
\begin{equation}\label{linearizedeigenvalueonbal01}
\left\{
\begin{array}{ll}
L_{\lambda}z=\tau z\,\, &\text{in}\,\, B_1^c,\\
z=0 &\text{on}\,\, \partial B_1,
\end{array}
\right.
\end{equation}
where $L_{\lambda}:=-\lambda\Delta +I-pu_\lambda^{p-1}$.
From Proposition \ref{prop21} we know that in $H_{0,r}^1\left(B_1^c\right)$ there is a unique negative eigenvalue $\tau_0$ to problem (\ref{linearizedeigenvalueonbal01}).
Let $z_\lambda$ be the positive eigenfunction corresponding to $\tau_0$ with $\left\Vert z_\lambda\right\Vert_{L^2\left(B_1^c\right)}=1$.

Define the functional $Q_\lambda:H_{0,G}^1\left(B_1^c\right)\rightarrow \mathbb{R}$ associated to
the above linearized eigenvalue problem
\begin{equation}
Q_\lambda\left(\psi\right)=\int_{B_1^c}\left(\lambda\vert\nabla \psi\vert^2+\psi^2-pu_\lambda^{p-1}\psi^2\right)\,\text{d}x.\nonumber
\end{equation}
Then we have that
\begin{proposition}\label{prop22} If $\lambda\leq-\tau_0S^2(1)/\mu_{i_1}$, there exists $\psi\in H_{0,G}^1\left(B_1^c\right)$ such that
$\int_{B_1^c} \psi z_\lambda\,\text{d}x=0$ and $Q_\lambda(\psi)\leq0$.
\end{proposition}
\textbf{Proof.} Let $\phi$ be an eigenfunction corresponding to $\mu_{i_1}$ with $\left\Vert \phi\right\Vert_{L^2\left(\mathbb{S}^{N-1}\right)}=1$.
Set $\psi:= z_\lambda(r)\phi(\theta)$. Then there holds that $\int_{B_1^c} \psi z_\lambda\,\text{d}x=0$. From the definition of $Q_\lambda$ we can obtain that
\begin{eqnarray}
Q_\lambda(\psi)&=&\int_{B_1^c}\left(\lambda\vert\nabla \psi\vert^2+\psi^2-pu_\lambda^{p-1}\psi^2\right)\,\text{d}x\nonumber\\
&=&\int_{B_1^c}\left(\lambda\vert\nabla z\vert^2\phi^2+\lambda\frac{1}{S^2}\vert\nabla \phi\vert^2z_\lambda^2+\psi^2-pu_\lambda^{p-1}\psi^2\right)\,\text{d}x\nonumber\\
&=&\int_{B_1^c}\left(\lambda\vert\nabla z_\lambda\vert^2+z_\lambda^2-pu_\lambda^{p-1}z_\lambda^2\right)\phi^2\,\text{d}x+\int_{B_1^c}\left(\lambda\frac{1}{S^2}\vert\nabla \phi\vert^2z_\lambda^2\right)\,\text{d}x\nonumber\\
&=&\int_{1}^{+\infty}S^{N-1}(r)\left(\lambda z_\lambda'^2+z_\lambda^2-pu_\lambda^{p-1}z_\lambda^2\right)\,\text{d}r\int_{\mathbb{S}^{N-1}}\phi^2\,\text{d}\theta\nonumber\\
& &+\lambda\int_{1}^{+\infty}S^{N-3}(r)z_\lambda^2\,\text{d}r\int_{\mathbb{S}^{N-1}}\vert\nabla \phi\vert^2\,\text{d}\theta\nonumber\\
&=&\tau_0+\lambda\int_{1}^{+\infty}S^{N-3}(r)z_\lambda^2\,dr\int_{\mathbb{S}^{N-1}}\vert\nabla \phi\vert^2\,\text{d}\theta\nonumber\\
&\leq&\tau_0+\frac{\mu_{i_1}\lambda}{S^2(1)},\nonumber
\end{eqnarray}
where $\mu_{i_1}=i_{1}(i_{1}+N-2)$.
If $\lambda\leq -\tau_0S^2(1)/\mu_{i_1}$, then we obtain that $Q_\lambda(\psi)\leq0$.\qed\\

Define
\begin{equation}
\Lambda_0=\sup\left\{\lambda>0:Q_\lambda\left(\psi\right)\leq0\,\,\text{for some}\,\,\psi\in H_{0,G}^1\left(B_1^c\right)\setminus\{0\}, \int_{B_1^c}\psi z_\lambda\,\text{d}x=0 \right\}.\nonumber
\end{equation}
It follows from Proposition \ref{prop22} that $\Lambda_0\geq -\tau_0S^2(1)/\mu_{i_1}>0$ and
Proposition \ref{prop22} implies that
$Q_\lambda$ can be negative if $\lambda$ is small.
The next conclusion shows that $Q_\lambda$ must be positive if $\lambda$ is large.
\begin{proposition}\label{prop23}
There exists $M >\Lambda_0$ such that for any $\lambda > M$, ${Q}_\lambda(\psi) > 0$ for any
$\psi\in E_0\setminus\{0\}$ where $E_0=\left\{\psi\in H_{0,G}^1\left(B_1^c\right):\int_{ B_1^c}\psi z_\lambda\,\text{d}x=0\right\}$.
\end{proposition}
We will prove this result as a corollary of a stronger result that we will need to use also later in our paper. More precisely, we consider the following more general functional on $H_G^1\left(B_1^c\right)$
\begin{equation}
\widetilde{Q}_\lambda\left(\psi\right)=\int_{B_1^c}\left(\lambda\vert\nabla \psi\vert^2+\psi^2-pu_\lambda^{p-1}\psi^2\right)\,\text{d}x-\lambda(N-1)\frac{e^2+1}{e^2-1}\int_{\partial B_1^c}\psi^2\,\text{d}s_x.\nonumber
\end{equation}
Observe that
\begin{equation}
\widetilde{Q}_\lambda\big|_{H_{0,G}^1\left(B_1^c\right)}={Q}_\lambda,\nonumber
\end{equation}
so Proposition \ref{prop23} is a corollary of the following result:
\begin{proposition}\label{prop24}
There exists $M >\Lambda_0$ such that for any $\lambda > M$, $\widetilde{Q}_\lambda(\psi) > 0$ for any
$\psi\in E\setminus\{0\}$ where $E=\left\{\psi\in H_{G}^1\left(B_1^c\right):\int_{\partial B_1}\psi\,\text{d}s_x=0, \int_{ B_1^c}\psi z_\lambda\,\text{d}x=0\right\}$.
\end{proposition}


The main objective of this section is to prove this Proposition \ref{prop24}. 
Although we adopt the idea of \cite[Proposition 5.3]{RRS}, there are some new
difficulties related to the hyperbolic space.

For $\psi\in H_G^1\left(B_R^c\right)$, define the functional
\begin{equation}
\widehat{Q}_R(\psi)=\int_{B_R^c}\left(\vert\nabla \psi\vert^2+\psi^2-pu_R^{p-1}\psi^2\right)\,\text{d}x-\frac{(N-1)(e^2+1)}{R(e^2-1)}\int_{\partial B_R^c}\psi^2\,\text{d}s_x . \nonumber
\end{equation}
According to the relations of $u_R$ and $u_\lambda$ with $\lambda=1/R^2$, we can see that
\begin{equation}
\widehat{Q}_R(\psi)=R^N\widetilde{Q}_\lambda(\psi).\nonumber
\end{equation}
So it is sufficient to show
$\widehat{Q}_R(\psi)$ is positive on $E_R$ for $R$ is sufficiently small, where
\begin{equation}
E_R=\left\{\psi\in H_{G}^1\left(B_R^c\right):\int_{\partial B_R}\psi\,\text{d}s_x=0,\int_{\partial B_R^c}\psi z_R\,\text{d}s_x=0\right\}\nonumber
\end{equation}
with $z_R(x)=z_\lambda\left(x/R\right)$. For this purpose, we first introduce some important results \cite{Ganguly} of the Schr\"{o}dinger equation on the whole space $\mathbb{H}^N$ as follows.
\begin{lemma}\label{prop31}
The following Schr\"{o}dinger equation
\begin{equation}
-\Delta u+u=u^p,\,\,u>0\,\,\text{in}\,\,\mathbb{H}^N,\nonumber
\end{equation}
possesses a unique radially symmetric $C^\infty$ solution $U$, which decreases and exponentially decays. Moreover, $U$ is non-degenerate in the space $H_{r}^1\left(\mathbb{H}^N\right)$ and has Morse index $1$.
\end{lemma}

It's worth noting that the nondegeneracy and radial Morse index one in Lemma \ref{prop31} are available in \cite[Theorem 3.1, Theorem 3.3]{Ganguly}.
\begin{lemma}\label{lemma31}
Let $u_n$ be the unique positive radial solution of (\ref{equationBR}) for $R = R_n$ with $\lim_{n\rightarrow+\infty}R_n=0$,
and $z_n =z_{R_n}/\left\Vert z_{R_n}\right\Vert$. Let us consider those functions extended to $
\mathbb{H}^N$ by 0. Then $u_n \rightarrow U$ and
$z_n \rightarrow Z$ in $H^1\left(\mathbb{H}^N\right)$, where $U$ is the radial ground state solution of the problem
\begin{equation}
-\Delta U+U=U^p,\,\,U>0\,\,\text{in}\,\,\mathbb{H}^N,\nonumber
\end{equation}
and $Z$ is the normalized positive eigenfunction corresponding to the negative eigenvalue
of the linearized problem.
\end{lemma}
\textbf{Proof.} Since the argument is the same as \cite{RRS} except that using the compact embedding results of \cite[Theorem 1.1]{Kazdan} or \cite[Lemma 4.2]{Berchio} and the uniqueness of \cite{Banica, Selvitella, Wang}, we omit it here.\qed
\\

Now we first give a useful inequality, which would be used later.
\begin{lemma}
For any $g\in C_0^{\infty}(\mathbb{R})$, $N \geq 2$ and $\lambda, r > 0$,
\begin{equation}
S^{N-2}(r) g^2(r)\leq\frac{1}{\lambda}\int_r^{+\infty}g'^2(s) S^{N-1}(s)\,\text{d}s+(2-N+\lambda)\int_r^{+\infty}g^2(s)S^{N-3}(s)\,\text{d}s,\nonumber
\end{equation}
where $S(s)=\sinh(s)$.
\end{lemma}
\textbf{Proof.} 
Observe that
\begin{eqnarray}
S^{N-2}(r) g^2(r)&=& (2-N)\int_r^{+\infty}S^{N-3}(s)S'(s)g^2(s)\,\text{d}s-2\int_r^{+\infty}S^{N-2}(s)g(s)g'(s)\,\text{d}s\nonumber\\
&\leq&(2-N)\int_r^{+\infty}S^{N-3}(s)g^2(s)\,\text{d}s-2\int_r^{+\infty}S^{N-2}(s)g(s)g'(s)\,\text{d}s.\nonumber
\end{eqnarray}
Using the Cauchy-Schwarz inequality as that of \cite[Lemma 7.2]{RRS} we can get that
\begin{eqnarray}
2\int_r^{+\infty}S^{N-2}(s)\left\vert g(s)g'(s)\right\vert\,\text{d}s&\leq& \lambda  \int_r^{+\infty}S^{N-3}(s)g^2(s)\,\text{d}s+\frac{1}{\lambda }\int_r^{+\infty}g'^2(s) S^{N-1}(s)\,\text{d}s,\nonumber
\end{eqnarray}
which implies the desired conclusion.\qed
\\

We will need the previous Lemma 
with $\lambda=\frac{4(N-1)}{3}$, that is
\[
\frac{4(N-1)}{3}\, S(r)^{N-2} g^2(r)\leq \int_r^{+\infty}g'^2(s) S^{N-1}(s)\,\text{d}s+\frac{4(N-1)(N+2)}{9}\int_r^{+\infty}g^2(s)S^{N-3}(s)\,\text{d}s\,.
\]

Using the above lemma we can get the following inequality.
\begin{lemma}\label{lemma33}
If the group of symmetries $G$ satisfies $(G1)$, then
\begin{eqnarray}
\frac{1}{S(R)}\int_{\partial B_R}\psi^2\,\text{d}s_x\leq \frac{3}{4(N-1)}\int_{B_R^c}\vert \nabla\psi\vert^2\,\text{d}x\nonumber
\end{eqnarray}
for any $\psi\in H_G^1\left(B_R^c\right)$ with $\int_{\partial B_R}\psi\,\text{d}s_x=0$ and $\int_{B_R^c}\psi z_R\,\text{d}x=0$.
\end{lemma}
\textbf{Proof.} By the Fourier expansion we have that
\begin{eqnarray}
\psi(r,\theta)=\sum_{k=1}^{+\infty}\psi_k(r)\phi_k(\theta),\nonumber
\end{eqnarray}
where $\phi_k$ are eigenfunctions of $-\Delta_{\mathbb{S}^{N-1}}$ with $G$-symmetry corresponding to the eigenvalues $\mu_{i_k}$.
Suppose that $G$ satisfies $(G1)$. Then
\begin{eqnarray}\label{lemma33inequality}
\int_{B_R^c}\left\vert \nabla\left(\psi_k(r)\phi_k(\theta)\right)\right\vert^2\,\text{d}r \text{d}\theta&=&\int_R^{+\infty}\left(\left(\psi_k'\right)^2S^{N-1}(r)+\mu_{i_k}\psi_k^2S^{N-3}(r)\right)\,\text{d}r
\int_{\partial B_1}\phi_k^2\,\text{d}\theta\nonumber\\
&\geq & \int_R^{+\infty}\left(\left(\psi_k'\right)^2S^{N-1}(r)+\mu_{i_1}\psi_k^2S^{N-3}(r)\right)\,\text{d}r
\int_{\partial B_1}\phi_k^2\,\text{d}\theta\nonumber\\
&\geq &\int_R^{+\infty}\left(\left(\psi_k'\right)^2S^{N-1}(r)+\frac{4(N+2)(N-1)}{9}\psi_k^2S^{N-3}(r)\right)\text{d}r\int_{\partial B_1}\phi_k^2\,\text{d}\theta\nonumber\\
&\geq &\frac{4(N-1)}{3}\, S(R)^{N-2} \psi_k^2(r)\, \int_{\partial B_1}\phi_k^2\,\text{d}\theta
\end{eqnarray}
where we have used the assumption $i_1 > \frac{2-N+\sqrt{(N-2)^2+\frac{16}{9}(N+2)(N-1)}}{2}$, which implies that
\[
\mu_{i_1}=i_1\left(i_1+N-2\right)\geq\frac{4(N+2)(N-1)}{9}\,,
\]
and for the last inequality we used the previous Lemma.
Note that
\begin{eqnarray}
\int_{\partial B_R}\left\vert \psi_k(r) \phi_k(\theta)\right\vert^2\,\text{d}s_x= S^{N-1}(R)\psi_k^2(R)\int_{\partial B_1}\phi_k^2(\theta)\,\text{d}\theta.\nonumber
\end{eqnarray}
This, combining with (\ref{lemma33inequality}), gives that
\begin{eqnarray}
\frac{1}{S(R)}\int_{\partial B_R}\left\vert \psi_k(r) \phi_k(\theta)\right\vert^2\,\text{d}s_x&\leq& \frac{3}{4(N-1)}\int_{B_R^c}\left\vert \nabla\left(\psi_k(r)\phi_k(\theta)\right)\right\vert^2\,\text{d}r \text{d}\theta,\nonumber
\end{eqnarray}
and the result is proved.\qed\\

Now we are ready to present the argument of the proof of Proposition \ref{prop24}.
\\
\\
\textbf{Proof of Proposition \ref{prop24}.} Take $R_n \to 0$, denote $B_n=B_{R_n}$, $u_n= u_{R_n}$ and $z_n= z_{R_n}$, and define:
$$ \chi_n = \inf \left \{ \widehat{Q}_{R_n}(\psi, \psi):\ \psi \in H^1_G(B_n^c),\ \int_{\partial B_n}\psi \text{d}s_x=0, \ \int_{B_n^c} \psi z_n\,\text{d}x =0,\ \int_{B_n^c} |\psi|^2\,\text{d}x= 1\right \}.$$
Assume, by contradiction, that $\chi_n \leq 0$. Arguing as \cite[Proposition 5.3]{RRS}, $\chi_n$ is attained.
Let
\begin{equation}
\chi_n=\widehat{Q}_{R_n}\left(\psi_n\right).\nonumber
\end{equation}
Then we have that $\psi_n$
is a solution of
\begin{equation}
-\Delta \psi+\psi-pu_n^{p-1}\psi=\chi_n\psi\,\, \text{in}\,\, B_n^c.
\nonumber
\end{equation}
We see that
\begin{equation}
\int_{B_n^c}\left(\left\vert\nabla \psi_n\right\vert^2+\psi_n^2-pu_n^{p-1}\psi_n^2\right)\,\text{d}x=\chi_n\int_{B_n^c}\psi_n^2\,\text{d}x+\frac{(N-1)\left(e^2+1\right)}{R_n\left(e^2-1\right)}\int_{\partial B_n}\psi_n^2\,\text{d}s_x.\nonumber
\end{equation}
Since $\psi_n$ is bounded, up to a subsequence, $\psi_n\rightharpoonup \psi_0\in H_G^1\left(B_r^c\right)$ for any $r > 0$, where $\psi_0\in H^1\left(\mathbb{H}^N\right)$.
By Lemma 3.1 and the argument of \cite[Proposition 5.3]{RRS}, we can derive that
\begin{equation}
\int_{B_n^c}pu_n^{p-1}\psi_n^2\,\text{d}x\rightarrow p\int_{\mathbb{H}^n}U^{p-1}\psi_0^2\,\text{d}x.\nonumber
\end{equation}
Now we estimate $\frac{1}{R_n}\int_{\partial B_n}\psi_n^2\,\text{d}s_x$.
Note that
\begin{equation}
\frac{1}{R_n}\int_{\partial B_n}\psi_n^2\,\text{d}s_x=\frac{1}{S\left(R_n\right)}\int_{\partial B_n}\psi_n^2\,\text{d}s_x+\left(\frac{1}{R_n}-\frac{1}{S\left(R_n\right)}\right)\int_{\partial B_n}\psi_n^2\,\text{d}s_x.\nonumber
\end{equation}
We have $\left( \frac{1}{R_n} - \frac{1}{S(R_n)}\right) \to 0$ when $R_n\to0$.
Since $\left\Vert \psi_n\right\Vert=1$ and the embedding $H^1\left(B_n^c\right)\hookrightarrow L^2\left(\partial B_n\right)$ is continuous (\cite{Kazdan} or \cite{Admas}), we see that
$\int_{\partial B_n}\psi_n^2\,\text{d}s_x$ is bounded.
Therefore, we have that
\begin{equation}
\lim_{n\rightarrow+\infty}\left(\frac{1}{R_n}-\frac{1}{S\left(R_n\right)}\right)\int_{\partial B_n}\psi_n^2\,\text{d}s_x=0.\nonumber
\end{equation}
So, for any $\varepsilon>0$, there exists an $N_0>0$ such that
\begin{equation}
\frac{(N-1)\left(e^2+1\right)}{ e^2-1}\left(\frac{1}{R_n}-\frac{1}{S\left(R_n\right)}\right)\int_{\partial B_n}\psi_n^2\,\text{d}s_x<\varepsilon\nonumber
\end{equation}
for any $n>N_0$.
It follows from Lemma \ref{lemma33} that if $G$ satisfies $(G1)$ then
\begin{eqnarray}
\frac{1}{S\left(R_n\right)}\int_{\partial B_n}\psi_n^2\,\text{d}s_x\leq \frac{3}{4(N-1)}\int_{B_n^c}\left\vert \nabla{\psi}_n(x)\right\vert^2\,\text{d}x\leq \frac{3}{4(N-1)}.\nonumber
\end{eqnarray}
So
\begin{equation}\label{boundedofka}
\int_{B_n^c}\left(\left\vert\nabla \psi_n\right\vert^2+\psi_n^2-pu_n^{p-1}\psi_n^2\right)\,\text{d}x-\chi_n\int_{B_n^c}\psi_n^2\,\text{d}x\leq\frac{3\left(e^2+1\right)}{4\left(e^2-1\right)}+\varepsilon
\end{equation}
for $n$ large enough.
It follows from $\chi_n\leq0$ that
\begin{eqnarray}
1-p\int_{B_n^c}u_n^{p-1}\psi_n^2\,\text{d}x\leq\frac{3\left(e^2+1\right)}{4\left(e^2-1\right)}+\varepsilon < 1.\nonumber
\end{eqnarray}
if $n$ is large enough. For the previous inequalities we used that $e^2\in(7,8)$.
As $R_n\rightarrow0$, we find that
\begin{eqnarray}
 p\int_{\mathbb{H}^n}U^{p-1}\psi_0^2\,\text{d}x>0,\nonumber
\end{eqnarray}
which implies that $\psi_0\not\equiv0$.
It follows from (\ref{boundedofka}) that
\begin{equation}
\chi_n\geq \frac{1- p\int_{B_n^c}u_n^{p-1}\psi_n^2\,\text{d}x-\frac{3\left(e^2+1\right)}{4\left(e^2-1\right)}-\varepsilon}{\int_{B_n^c}\psi_n^2\,\text{d}x}.\nonumber
\end{equation}
Since
\begin{equation}
\liminf_{n\rightarrow+\infty}\int_{B_n^c}\psi_n^2\,\text{d}x\geq \int_{\mathbb{H}^n}\psi_0^2\,\text{d}x>0,\nonumber
\end{equation}
it's obvious that
\begin{equation}
\chi_n\geq -\frac{p\int_{\mathbb{H}^n}U^{p-1}\psi_0^2\,\text{d}x}{\int_{\mathbb{H}^n}\psi_0^2\,\text{d}x}.\nonumber
\end{equation}
Thus $\chi_n$ is bounded and there exists $\chi_0\leq0$ such that $\lim_{n\rightarrow+\infty}\chi_n=\chi_0$.
Observe that $\psi_n$ is the weak solution
\begin{equation}
-\Delta \psi+\psi-pu_n^{p-1}\psi=\chi_n\psi\,\, \text{in}\,\, B_n^c.
\nonumber
\end{equation}
So $\psi_0$ is a nontrivial weak
solution of
\begin{equation}
-\Delta \psi+\psi-pU^{p-1}\psi=\chi_0 \psi\,\,\text{in}\,\,\mathbb{H}^N\setminus\{0\}.\nonumber
\end{equation}
Since $\psi_0\in H^1\left(\mathbb{H}^N\right)$, the singularity is removable. So, it is a weak solution in the
whole $\mathbb{H}^N$.
Hence, we get that
\begin{equation}
\left\Vert \psi_0\right\Vert^2=p\int_{\mathbb{H}^N}U^{p-1}\psi_0^2\,\text{d}x+\int_{\mathbb{H}^N}\chi_0 \psi_0^2\,\text{d}x.\nonumber
\end{equation}
Let $\widetilde{\psi}_0(r)=\int_{\mathbb{S}^{N-1}}\psi_0\,d\theta$.
Then we have that
\begin{equation}
-\Delta \widetilde{\psi}_0+\widetilde{\psi}_0-pU^{p-1}\widetilde{\psi}_0=\chi_0 \widetilde{\psi}_0\,\,\text{in}\,\,(0,+\infty).\nonumber
\end{equation}
The only possibility is $\widetilde{\psi}_0= kZ$, $k\neq0$.
By the same arguments as in Step 3 of \cite[Proposition 5.3]{RRS}, we conclude that
\begin{equation}
0=\int_{B_n^c}\psi_nz_n\,\text{d}x\rightarrow \int_{\mathbb{H}^N}\psi_0Z\,\text{d}x=\int_0^{+\infty}S^{N-1}\widetilde{\psi}_0Z\,\text{d}r=k\int_{\mathbb{H}^N}Z^2\,\text{d}x\neq0.\nonumber
\end{equation}
which yields the desired contradiction.\qed\\

\textbf{Remark.} From the above argument we understand that the condition of $i_1$ of hypothesis $(G1)$ is to guarantee the inequality in the form of Lemma \ref{lemma33}, that is to show
\begin{eqnarray}
\frac{3\left(e^2+1\right)}{4\left(e^2-1\right)}<1\nonumber
\end{eqnarray}
in the previous proof. This is a difference with respect to the Euclidean space and is due to the nature of the hyperbolic space.\\

In addition, the behavior of the functional $\widetilde{Q}_\lambda(\psi)$ near $\Lambda_0$ can be also established as follows.
\begin{proposition}\label{prop25}
$\widetilde{Q}_{\Lambda_0}\left(\psi\right)<0$ for some $\psi\in E$.
\end{proposition}
\noindent\textbf{Proof.} Since the proof is the same as that of \cite[Lemma 5.4]{RRS} with obvious changes, we omit it here.\qed
\\

Combining Proposition \ref{prop24} with Proposition \ref{prop25}, we have that $\widetilde{Q}_\lambda$ can be negative for $\lambda$ near $\Lambda_0$, but is positive if $\lambda$ is large enough. This provides a possibility for finding the zero point of $\widetilde{Q}_\lambda$ with respect to $\lambda$.
We end this section by showing an existence result which will be used later.
\begin{lemma}\label{lemma25}
For $\lambda>\Lambda_0$, $L_\lambda$ is an isomorphism.
Given $v\in  H^{1/2}_G \left(\mathbb{S}^{N-1}\right)$, there exists a unique solution $\psi_v\in H_G^1\left(B_1^c\right)$ of the problem
\begin{equation}\label{eigenvalueonc1cvnu=0}
\left\{
\begin{array}{ll}
-\lambda\Delta\psi+\psi-pu_\lambda^{p-1}\psi=0\,\, &\text{in}\,\, B_1^c,\\
\psi=v &\text{on}\,\, \partial B_1.
\end{array}
\right.
\end{equation}
In addition, $\int_{B_1^c}\psi_v z_\lambda\,\text{d}x=0$ and $\int_{\partial B_1}\frac{\psi_v}{\partial \nu} \,\text{d}s_x=0$.
\end{lemma}
\noindent\textbf{Proof.} Since the proof is the same as that of \cite[Lemma 3.5, Lemma 4.2]{RRS} with obvious changes, we omit it here.\qed
\\ \\
\textbf{Remark.} Using the spherical coordinates $(r,\omega)\in [0,+\infty)\times \mathbb{S}^{N-1}$, the equation (\ref{eigenvalueonc1cvnu=0}) can be written as
\begin{equation}
\left\{
\begin{array}{ll}
-\lambda\left(\partial_r^2+(N-1)\frac{\cosh r}{\sinh r}\partial_r+\frac{1}{\sinh^2 r}\Delta_{\mathbb{S}^{N-1}}\right)\psi+\psi-pu_\lambda^{p-1}\psi=0\,\, &\text{in}\,\, \Omega,\\
\psi=v &\text{on}\,\, \partial \Omega,
\end{array}
\right.\nonumber
\end{equation}
where $\Omega=(0,+\infty)\times \mathbb{S}^{N-1}\subset \mathbb{R}^N$. Hence the standard Schauder elliptic estimates \cite{Gilbarg}
is valid for (\ref{eigenvalueonc1cvnu=0}). In particular, if $v\in C^{2,\alpha}_{G,m}\left(\partial B_1\right)$, one has $\psi_v\in C^{2,\alpha}_{G}\left(B_1^c\right)$.

\section{Normal derivative operator and its linearization.}

For $\lambda>\Lambda_0$ and each $v\in C_{G,m}^{2,\alpha}\left(\mathbb{S}^{N-1}\right)$ whose norm is sufficient small,
we consider the following perturbation problem
\begin{equation}\label{operator}
\left\{
\begin{array}{ll}
-\lambda\Delta u+u-u^p=0\,\, &\text{in}\,\, B_{1+v}^c,\\
u=0 &\text{on}\,\, \partial B_{1+v}^c.
\end{array}
\right.
\end{equation}
By the nondegeneracy of $L_\lambda$, using an argument similar to that of \cite[Proposition 4.1]{RRS}, we have the following existence of a solution to problem (\ref{operator}).
\begin{proposition}\label{prop41}
For any $\lambda>\Lambda_0$ and all $v\in C_{G,m}^{2,\alpha}\left(\mathbb{S}^{N-1}\right)$ whose norm is sufficiently small, there exists a unique positive solution $u =u(\lambda,v)\in \mathcal{C}^{2,\alpha}\left({B_{1+v}^c}\right)\cap H_{0,G}^1\left({B_{1+v}^c}\right)$ to (\ref{operator}). Moreover, $u$ depends smoothly on $v$, $\lambda$ (in the sense that considering a regular diffeomorphism mapping $B_1^c$ into $B_{1+v}^c$ and considering the pullback of $u$ on $B_1^c$, then the map $(\lambda,v) \to u$ is smooth). In particular, $u = u_\lambda$ when $v \equiv0$.
\end{proposition}

We introduce the operator $F:\left(\Lambda_0,+\infty\right)\times C_{G,m}^{2,\alpha}\left(\mathbb{S}^{N-1}\right)\rightarrow C_{G,m}^{1,\alpha}\left(\mathbb{S}^{N-1}\right)$ defined by
\begin{equation}
F(\lambda,v)=\frac{1}{\text{Vol}\left(\partial B_{1+v}\right)}\int_{\partial B_{1+v}}\frac{\partial u(\lambda,v)}{\partial \nu}\,\text{dvol}-\frac{\partial u(\lambda,v)}{\partial \nu},\nonumber
\end{equation}
where $\nu$ denotes the unit normal vector field to $\partial B_{1+v}$ pointing to the interior of $B_{1+v}$.
Since $\partial_r u_\lambda(1)$ is a constant, then $F(\lambda,0)=0$ for any $\lambda>\Lambda_0$.
Therefore, finding the nontrivial domains emanating from $B_1^c\subset \mathbb{H}^N$ such that problem (\ref{ffunction}) has a bounded solution is equivalent to study the nontrivial   solutions of $F(\lambda,v)=0$ bifurcating from $F(\lambda,0)=0$.
It follows from the properties of $u(\lambda,v)$ that $F$ is $C^1$.

We introduce the operator
\begin{equation}\label{ksaifenl}
{H}_\lambda(v)={\left.\frac{\partial \psi_v}{\partial \nu}\right|_{\partial B_1}-(N-1)\frac{e^2+1}{e^2-1}v,}
\end{equation}
where $\psi_v$ is given by Lemma \ref{lemma25}.
As in \cite[Proposition 4.3]{RRS} we can show that the linearization of the operator $F$ with respect to $v$ at point $(\lambda,0)$ is $C{H}_\lambda$ where $C=u_\lambda'(1)\neq0$.
The only difference is that here
\[
u_\lambda''(1)=-(N-1) \frac{C(1)}{S(1)}u_\lambda'(1)\,.
\]
By Schauder elliptic estimates \cite[Theorem 8.16, Theorem 8.13, Theorem 6.6]{Gilbarg}, we can obtain that $\sup_{B_{1}^{c}}\left\vert \psi_v\right\vert\leq \sup_{\partial B_{1}^{c}}\vert v\vert$ and $\left\Vert \psi\right\Vert_{C_G^{2,\alpha}\left(B_{1}^{c}\right)}\leq  M \Vert v\Vert_{C_{G,m}^{2,\alpha}\left(\mathbb{S}^{N-1}\right)}$ for some $M>0$.
It follows that
\[
\left\Vert {H}_\lambda(v)\right\Vert_{C_{G,m}^{1,\alpha}\left(\mathbb{S}^{N-1}\right)}\leq M_1 \Vert v\Vert_{C_{G,m}^{2,\alpha}\left(\mathbb{S}^{N-1}\right)}
\]
for some $M_1>0$.
So, the operator
\begin{equation}
{H}_\lambda:C_{G,m}^{2,\alpha}\left(\mathbb{S}^{N-1}\right)\longrightarrow C_{G,m}^{1,\alpha}\left(\mathbb{S}^{N-1}\right)\nonumber
\end{equation}
is bounded.
Furthermore, we also have the:
\begin{proposition}\label{prop42}
The operator ${H}_\lambda:C_{G,m}^{2,\alpha}\left(\mathbb{S}^{N-1}\right)\longrightarrow C_{G,m}^{1,\alpha}\left(\mathbb{S}^{N-1}\right)$ is a self-adjoint, first order elliptic operator.
\end{proposition}
\textbf{Proof.} Since $H_\lambda$ is the sum of the Dirichlet-to-Neumann operator for $-\lambda\Delta+1-pu_\lambda^{p-1}$ and
a constant times the identity, ${H}_\lambda$ is a first order elliptic operator.
So it is enough to prove it is self-adjoint. Let $\psi_1$ and $\psi_2$ be the solution of problem (\ref{eigenvalueonc1cvnu=0})
with $v=v_1$ and $v=v_2$.
Then we have that
\begin{eqnarray}
\int_{\mathbb{S}^{N-1}}\left(H_\lambda\left(v_1\right)v_2-H_\lambda\left(v_2\right)v_1\right)\,\text{d}\theta
=\int_{\mathbb{S}^{N-1}}\left(\frac{\partial \psi_1}{\partial \nu}v_2-\frac{\partial \psi_2}{\partial \nu}v_1\right)\,\text{d}\theta.\nonumber
\end{eqnarray}
We multiply the equation of $\psi_1$ by $\psi_2$ and the equation of $\psi_2$ by $\psi_1$, integrating by parts, we obtain that
\begin{eqnarray}
\int_{\mathbb{S}^{N-1}}\left(\frac{\partial \psi_1}{\partial \nu}v_2-\frac{\partial \psi_2}{\partial \nu}v_1\right)\,\text{d}\theta=0.\nonumber
\end{eqnarray}
So \begin{eqnarray}
\int_{\mathbb{S}^{N-1}}\left(H_\lambda\left(v_1\right)v_2-H_\lambda\left(v_2\right)v_1\right)\,\text{d}\theta=0,\nonumber
\end{eqnarray}
which verifies $H_\lambda$ is self-adjoint.\qed
\\

The necessary condition to bifurcate is that ${H}_\lambda$ degenerates.
So we will find the value of $\lambda$ such that ${H}_\lambda$ is degenerate.
For any $v\in C_{G,m}^{2,\alpha}\left(\mathbb{S}^{N-1}\right)$, by virtue of the Fourier expansion with respect to spherical harmonics \cite[Theorem 3.2.11]{Groemer}, $v$ can be written as
\begin{equation}
v=\sum_{k=1}^\infty\sum_{j=1}^{m_k}a_{i_k,j}\zeta_{i_k,j}(\theta),\nonumber
\end{equation}
where $\zeta_{i_k,j}$ (normalized to $1$ in the $L^2$-norm) is an eigenfunction corresponding to $\mu_{i_k}$, being
\[
\text{span}\left\{\zeta_{i_k,1}, \ldots,\zeta_{i_k,m_k}\right\}
\]
the associate eigenspace. We need now to study the eigenvalues of ${H}_\lambda$.
\begin{proposition} \label{prop43}
For any $\lambda>\Lambda_0$, ${H}_\lambda$ possesses a sequence of eigenvalues $\left\{\sigma_{i_k}(\lambda)\right\}_{k\in \mathbb{N}}$ such that
\begin{equation}
\sigma_{i_1}(\lambda)<\cdots<\sigma_{i_k}(\lambda)<\cdots.\nonumber
\end{equation}
The eigenfunction corresponding to $\sigma_{i_k}(\lambda)$ is $\sum_{j=1}^{m_k}a_{i_k,j}\zeta_{i_k,j}(\theta)$ with $\sum_{j=1}^{m_k}a_{i_k,j}^2\neq0$.
\end{proposition}
\textbf{Proof.} Let $\phi_0(r,\theta)=u_\lambda'(r)v\left(\theta\right)$ with $r=\vert x\vert$.
We can verify that
\begin{equation}
\Delta u_\lambda'=(N-1)\left(\frac{S'^2}{S^2}-1\right)u_\lambda'+\frac{u_\lambda'-pu_\lambda^{p-1}u_\lambda'}{\lambda}.\nonumber
\end{equation}
So we have that
\begin{eqnarray}
\Delta \phi_0&=&v\Delta u_\lambda'+u_\lambda'\Delta v\nonumber\\
&=&\left((N-1)\left(\frac{S'^2}{S^2}-1\right)u_\lambda'+\frac{u_\lambda'-pu_\lambda^{p-1}u_\lambda'}{\lambda}\right)v+
u_\lambda'\Delta v\nonumber\\
&=&\sum_{k=1}^\infty\sum_{j=1}^{m_k}a_{i_k,j}\zeta_{i_k,j}(\theta)\left((N-1)\left(\frac{S'^2}{S^2}-1\right)u_\lambda'
+\frac{u_\lambda'-pu_\lambda^{p-1}u_\lambda'}{\lambda}\right)\nonumber\\
& &-\frac{1}{S^2}\sum_{k=1}^\infty\sum_{j=1}^{m_k}a_{i_k,j}\mu_{i_k}\zeta_{i_k,j}(\theta)
u_\lambda'\nonumber\\
&=&\sum_{k=1}^\infty\sum_{j=1}^{m_k}a_{i_k,j}u_\lambda'\zeta_{i_k,j}(\theta)\left((N-1)\left(\frac{S'^2}{S^2}-1\right)
+\frac{1-pu_\lambda^{p-1}}{\lambda}-\frac{\mu_{i_k}}{S^2}\right).\nonumber
\end{eqnarray}
It follows that
\begin{eqnarray}
-\lambda \Delta \phi_0+\phi_0-pu_\lambda^{p-1}\phi_0=\sum_{k=1}^\infty\sum_{j=1}^{m_k}a_{i_k,j}u_\lambda'\zeta_{i_k,j}(\theta)\lambda\left(\frac{\mu_{i_k}}{S^2}-(N-1)\left(\frac{S'^2}{S^2}-1\right)
\right).\nonumber
\end{eqnarray}

Consider the following problem
\begin{eqnarray}
\left\{
\begin{array}{lll}
-\lambda \Delta \Psi+\Psi-pu_\lambda^{p-1}\Psi=\displaystyle  \sum_{k=1}^\infty\sum_{j=1}^{m_k}a_{i_k,j}u_\lambda'\zeta_{i_k,j}(\theta)\lambda\left(\frac{\mu_{i_k}}{S^2}-(N-1)\left(\frac{S'^2}{S^2}-1\right)
\right)\,\, &\text{in}\,\, B_1^c,\\
\Psi=0  &\text{on}\,\, \partial B_1.
\end{array}
\right.\nonumber
\end{eqnarray}
By the Riesz Theorem, the above problem has a unique solution $\Psi$ in $H^1(B_1^c)$.
Then we see that
\begin{equation}\label{phixt}
\Psi(r,\theta)= {\phi_0(r,\theta)}-\psi_v u_\lambda'(1).
\end{equation}
Moreover, we have that
\begin{equation}
\partial_r \Psi(1,\theta)= {u_\lambda''(1)v}-u_\lambda'(1) \partial_r\psi_v(1,\theta).\nonumber
\end{equation}
Note that $u_\lambda''(1)=-(N-1)u_\lambda'(1)C(1)/S(1)$.
We have that
\begin{equation}
\partial_r \Psi(1,\theta)={-(N-1)\frac{e^2+1}{e^2-1}u_\lambda'(1)v}-u_\lambda'(1) \partial_r\psi_v(1,\theta)={u_\lambda'(1) \left[ \left.\partial_\nu\psi_v\right|_{\partial B_1} -(N-1)\frac{e^2+1}{e^2-1}v \right] = u_\lambda'(1)H_\lambda v}\nonumber
\end{equation}
{because on $\partial B_1$ we have $\partial_\nu = -\partial_r$.}
Let $V_k$ be the space spanned by the functions $\zeta_{i_k,1}(\theta)$, $\ldots$, $\zeta_{i_k,m_k}(\theta)$. Since
\[
H_\lambda v = \frac{\partial_r \Psi(1,\theta)}{u_\lambda'(1)}\,,
\]
we deduce that $H_\lambda$ preserves $V_k$ for any $k$.
It follows that
\begin{equation}
{H}_\lambda \left(\sum_{j=1}^{m_k}a_{i_k,j}\zeta_{i_k,j}(\theta) \right)=\sigma_{i_k}(\lambda)\sum_{j=1}^{m_k}a_{i_k,j}\zeta_{i_k,j}(\theta),\nonumber
\end{equation}
where $\sigma_{i_k}(\lambda)$ are the eigenvalues of $H_\lambda$, and $\sum_{j=1}^{m_k}a_{i_k,j}\zeta_{i_k,j}(\theta)$ with $\sum_{j=1}^{m_k}a_{i_k,j}^2\neq0$ are the eigenfunctions associated to $\sigma_{i_k}(\lambda)$.
So we have that
\begin{equation}
{H}_\lambda \left( \sum_{k=1}^n\sum_{j=1}^{m_k}a_{i_k,j}\zeta_{i_k,j}(\theta)\right)=\sum_{k=1}^n\sigma_{i_k}(\lambda)\sum_{j=1}^{m_k}a_{i_k,j}\zeta_{i_k,j}(\theta).\nonumber
\end{equation}
Since $H_\lambda$ is bounded, we have that
\begin{eqnarray}\label{ksaifenli1}
{H}_\lambda (v)&=&\lim_{n\rightarrow+\infty}{H}_\lambda \left( \sum_{k=1}^{n}\sum_{j=1}^{m_k}a_{i_k,j}\zeta_{i_k,j}(\theta)\right)=\sum_{k=1}^{+\infty}\sigma_{i_k}(\lambda)\sum_{j=1}^{m_k}a_{i_k,j}\zeta_{i_k,j}(\theta).
\end{eqnarray}
Next we study the properties of $\sigma_{i_k}(\lambda)$.
With the help of (\ref{phixt}), we can write $\Psi(r,\theta)$ as
\begin{equation}
\Psi(r,\theta)=\sum_{k=1}^\infty\sum_{j=1}^{m_k}b_{i_k}(r)a_{i_k,j}\zeta_{i_k,j}(\theta),\nonumber
\end{equation}
where $b_{k}$ is the continuous solution on $[1,+\infty)$ of
\begin{equation}
-\lambda\left(\partial_r^2+(N-1)\frac{C(r)}{S(r)}\partial_r\right)b +\lambda \frac{\mu_{i_k}}{S^2(r)}b+\left(1-pu_\lambda^{p-1}\right)b=
\lambda u_\lambda'\left(\frac{\mu_{i_k}}{S^2}-(N-1)\left(\frac{S'^2}{S^2}-1\right)\right)\nonumber
\end{equation}
with $b_{i_k}(1)=0$ and $\lim_{r\rightarrow+\infty}b_{i_k}(r)=0$.
Hence, we conclude that
\begin{equation}
u_\lambda'(1)\sigma_{i_k}(\lambda)=\partial_r b_{i_k}(1)\nonumber
\end{equation}
for any $k\in \mathbb{N}$.
Combining (\ref{eigenvalueonc1cvnu=0}), (\ref{ksaifenl}) with (\ref{ksaifenli1}) we derive that
\begin{equation}
\psi_v=\sum_{{k=1}}^\infty\sum_{j=1}^{m_k}c_{i_k}(r)a_{i_k,j}\zeta_{i_k,j}(\theta),\nonumber
\end{equation}
where $c_{i_k}$ is the continuous solution on $[1,+\infty)$ of
\begin{equation}\label{ckequation}
-\lambda\left(\partial_r^2+(N-1)\frac{C(r)}{S(r)}\partial_r\right)c+\lambda\frac{\mu_{i_k}}{S^2(r)}c +\left(1-pu_\lambda^{p-1}\right)c=0
\end{equation}
with $c_{i_k}(1)=1$ and $\lim_{r\rightarrow+\infty}c_{i_k}(r)=0$.
This implies that
\begin{equation}
\sigma_{i_k}(\lambda)=-\left(c_{i_k}'(1)+(N-1)\frac{e^2+1}{e^2-1}\right)\nonumber
\end{equation}
for $k\geq1$. 

When $k=1$, let $\widetilde{c}(r,\theta)=c_{i_1}(r)\phi(\theta)$ with $\theta\in \mathbb{S}^{N-1}$, where
$\phi$ is a normalized eigenfunction corresponding to $\mu_{i_1}$.
So, using (\ref{ckequation}), we have that
\begin{equation}
\Delta \widetilde{c}=S^{1-N}(r) \left(S^{N-1}(r)c'_{i_1}\right)'\phi(\theta)+\frac{c_{i_1}}{S^2}\Delta \phi(\theta)=\frac{1-pu_\lambda^{p-1}}{\lambda}\widetilde{c},\nonumber
\end{equation}
where the first equal sign is obtained by applying \cite[Formala (31) in Chapter II]{Chavel}.
Multiplying $\widetilde{c}$ on both sides and integrating it by part, we get that
\begin{equation}
-\int_{\partial B_1}c\nabla \widetilde{c}\cdot\nu\,\text{d}s_x+\int_{B_1^c} \left\vert\nabla \widetilde{c}\right\vert^2\,\text{d}x+\int_{B_1^c}\frac{1-pu_\lambda^{p-1}}{\lambda}\widetilde{c}^2\,\text{d}x=0.\nonumber
\end{equation}
It follows that
\begin{equation}
-c_{i_1}'(1)=\int_1^\infty S^{N-1}\left(c_{i_1}'^2+\frac{1-pu_\lambda^{p-1}}{\lambda}c_{i_1}^2\right)\,\text{d}r+\mu_{i_1}\int_1^\infty S^{N-3}c_{i_1}^2\,\text{d}r.\nonumber
\end{equation}
By definition, it's known that
\begin{equation}
\widetilde{Q}_\lambda\left(\psi\right)={Q}_\lambda\left(\psi\right)-\lambda(N-1)\frac{e^2+1}{e^2-1}\int_{\partial B_1}\psi^2\,\text{d}s_x\nonumber
\end{equation}
and we define
\begin{equation}
E=\left\{\phi\in H_G^1\left(B_1^c\right):\int_{\partial B_1}\phi\,\text{d}s_x=0,\int_{B_1^c}\phi z_\lambda\,\text{d}x=0\right\}\nonumber
\end{equation}
with \begin{equation}
H_{G}^1\left(B_1^c\right)=\left\{u\in H^1\left(B_1^c\right):u=u\circ g, \forall g\in G\right\}.\nonumber
\end{equation}
For any $\psi\in E$ with $\int_{\partial B_1}\psi^2\,\text{d}s_x=1$,
there exist $\psi_k$ ({$k\geq 1$}) and $a_{i_k,j}$ ($j\in\left\{1,\ldots,m_k\right\}$) with $\sum_{j=1}^{m_k}a_{i_k,j}^2=1$ for any $k$ such that
\begin{equation}
\psi(r,\theta)=\sum_{{k=1}}^{+\infty}\psi_k(r)\sum_{j=1}^{m_k}a_{i_k,j}\zeta_{i_k,j}(\theta).\nonumber
\end{equation}
Then we have that
\begin{eqnarray}
\frac{1}{\lambda}\widetilde{Q}_\lambda \left(\psi\right)&=&\sum_{k=1}^{+\infty}\int_1^\infty S^{N-1}\left(\psi_{i_k}'^2+\frac{1-pu_\lambda^{p-1}}{\lambda}\psi_{i_k}^2\right)\,\text{d}r+\sum_{k=1}^{+\infty}\mu_{i_k}\int_1^\infty S^{N-3}\psi_{i_k}^2\,\text{d}r\nonumber\\
& &-(N-1)\frac{e^2+1}{e^2-1}\nonumber\\
&\geq&\sum_{k=1}^{+\infty}\int_1^\infty S^{N-1}\left(\psi_{i_k}'^2+\frac{1-pu_\lambda^{p-1}}{\lambda}\psi_{i_k}^2\right)\,\text{d}r+\mu_{i_1}\sum_{k=1}^{+\infty}\int_1^\infty S^{N-3}\psi_{i_k}^2\,\text{d}r\nonumber\\
& &-(N-1)\frac{e^2+1}{e^2-1}
\nonumber\\
&\geq&\frac{1}{\lambda}\widetilde{Q}_\lambda\left(\widetilde{\phi}\right),\nonumber
\end{eqnarray}
where
\begin{equation}
\widetilde{\phi}(r,\theta)=\sum_{k=1}^{+\infty}\psi_k(r)\sum_{j=1}^{m_1}a_{i_1,j}\zeta_{i_1,j}(\theta)\nonumber
\end{equation}
and $\int_{\partial B_1}\widetilde{\phi}^2\,\text{d}s_x=1$.
Therefore, the infimum of $\widetilde{Q}_\lambda(\phi)/\lambda$ in $E$ with $\int_{\partial B_1}\phi^2\,\text{d}s_x=1$ is attained.
Define
\begin{equation}
\sigma_1\left(H_\lambda\right):=\inf\left\{\frac{1}{\lambda}\widetilde{Q}_\lambda(\phi):\phi\in E,\int_{\partial B_1}\phi^2\,\text{d}s_x=1\right\}.\nonumber
\end{equation}
We next to investigate the relations of $\sigma_1\left(H_\lambda\right)$ and $\sigma_{i_1}(\lambda)$.

The above argument implies that there exists $\phi(r,\theta)=\upsilon(r)\sum_{j=1}^{m_1}a_{i_1,j}\zeta_{i_1,j}(\theta)$ with $\sum_{j=1}^{m_1}a_{i_1,j}^2=1$ such that
$\phi\in E$, $\int_{\partial B_1}\phi^2\,\text{d}s_x=1$ and
\begin{equation}
\sigma_1\left(H_\lambda\right)=\frac{1}{\lambda}\widetilde{Q}_\lambda(\phi).\nonumber
\end{equation}
It follows that
\begin{equation}
\sigma_1\left(H_\lambda\right)=\int_1^\infty S^{N-1}\left(\upsilon'^2+\frac{1-pu_\lambda^{p-1}}{\lambda}\upsilon^2\right)\,\text{d}r+\mu_{i_1}\int_1^\infty S^{N-3}\upsilon^2\,\text{d}r-(N-1)\frac{e^2+1}{e^2-1},\nonumber
\end{equation}
which is the functional of
\begin{equation}\label{weaksolution}
\left\{
\begin{array}{lll}
-\lambda\left(\partial_r^2+\frac{(N-1)C(r)}{S(r)}\partial_r\right)\upsilon+\lambda\frac{\mu_{i_1}}{S^2(r)}\upsilon+\left(1-pu_\lambda^{p-1}\right)\upsilon=0,\\
\upsilon(1)=1.
\end{array}
\right.
\end{equation}
So $\upsilon$ is a weak solution of (\ref{weaksolution}).
By Schauder elliptic estimates, $\upsilon$ is also the classical solution of (\ref{weaksolution}).
By Lemma 2.5 we deduce that $\upsilon(r)\equiv c_{i_1}(r)$.
So we get that
\begin{equation}
\sigma_1\left(H_\lambda\right)=-c_{i_1}'(1)-(N-1)\frac{e^2+1}{e^2-1}=\sigma_{i_1}(\lambda).\nonumber\nonumber
\end{equation}
Therefore, we obtain that $\sigma_{i_1}(\lambda)$ is the same as $\sigma_1\left(H_\lambda\right)$. So the eigenspace corresponding to $\sigma_{i_1}(\lambda)$
is just $V_1$.
We next study the high eigenvalues.

For $\psi\in E$, we call $\psi\in E_k$ if there exist $\psi_{i_k}$ and $a_{i_k,j}$ such that
\begin{equation}
\psi(r,\theta)=\psi_{i_k}(r)\sum_{j=1}^{m_k}a_{i_k,j}\zeta_{i_k,j}(\theta).\nonumber
\end{equation}
Set
\begin{equation}
E_{k-1}^c:=E\setminus \cup_{i=1}^{k-1}E_{i}.\nonumber
\end{equation}
We will explore the infimum of $\widetilde{Q}_\lambda(\phi)/\lambda$ in $E_{k-1}^c$ with $\int_{\partial B_1}\phi^2\,\text{d}S=1$.
For a function $\phi$ on $[1,+\infty)$ we define
\begin{equation}
Q_{\lambda,k}(\phi)=\int_1^\infty S^{N-1}\left(\phi'^2+\frac{1-pu_\lambda^{p-1}}{\lambda}\phi^2\right)\,\text{d}r+\mu_{i_k}\int_1^\infty S^{N-3}\phi^2\,\text{d}r-(N-1)\frac{e^2+1}{e^2-1}.\nonumber
\end{equation}
Since $\mu_{i_k}$ is increasing with respect to $k$, $Q_{\lambda,k}(\phi)$ is increasing with respect to $k$.
For any $\psi\in E_{{k-1}}^c$ with $\int_{\partial B_1}\psi^2\,\text{d}s_x=1$,
there exist $\psi_l$ ($l\geq k$) and $a_{i_l,j}$ ($j\in\left\{1,\ldots,m_l\right\}$) with $\sum_{j=1}^{m_l}a_{i_l,j}^2=1$ for any $l\geq k$ such that
\begin{equation}
\psi(r,\theta)=\sum_{l=k}^{+\infty}\psi_l(r)\sum_{j=1}^{m_l}a_{i_l,j}\zeta_{i_l,j}(\theta).\nonumber
\end{equation}
Then we have that
\begin{eqnarray}
\frac{1}{\lambda}\widetilde{Q}_\lambda\left(\psi\right)&=&\sum_{l=k}^{+\infty}\int_1^\infty S^{N-1}\left(\psi_l'^2+\frac{1-pu_\lambda^{p-1}}{\lambda}\psi_l^2\right)\,\text{d}r+\sum_{l=k}^{+\infty}\mu_{i_l}\int_1^\infty S^{N-3}\psi_{i_l}^2\,\text{d}r\nonumber\\
& &-(N-1)\frac{e^2+1}{e^2-1}\nonumber\\
&\geq&\sum_{l=k}^{+\infty}\int_1^\infty S^{N-1}\left(\psi_l'^2+\frac{1-pu_\lambda^{p-1}}{\lambda}\psi_l^2\right)\,\text{d}r+\mu_{i_k}\sum_{i=l}^{+\infty}\int_1^\infty S^{N-3}\psi_l^2\,\text{d}r\nonumber\\
& &-(N-1)\frac{e^2+1}{e^2-1}
\nonumber\\
&\geq&\frac{1}{\lambda}\widetilde{Q}_\lambda\left(\widetilde{\phi}\right),\nonumber
\end{eqnarray}
where
\begin{equation}
\widetilde{\phi}(r,\theta)=\sum_{l=k}^{+\infty}\psi_l(r)\sum_{j=1}^{m_k}a_{i_k,j}\zeta_{i_k,j}(\theta)\in E_k\nonumber
\end{equation}
and $\int_{\partial B_1}\widetilde{\phi}^2\,\text{d}s_x=1$.
Therefore, the infimum is attained in $E_k$.
Hence we can set
\begin{equation}
\sigma_k\left(H_\lambda\right):=\inf\left\{\frac{1}{\lambda}\widetilde{Q}_\lambda(\phi):\phi\in E_{k-1}^c,\int_{\partial B_1}\phi^2\,\text{d}s_x=1\right\}.\nonumber
\end{equation}

There exists $\phi(r,\theta)=\upsilon(r)\sum_{j=1}^{m_k}a_{i_k,j}\zeta_{i_k,j}(\theta)$ with $\sum_{j=1}^{m_k}a_{i_k,j}^2=1$ such that
$\phi\in E_{k-1}^c$, $\int_{\partial B_1}\phi^2\,\text{d}s_x=1$ and
\begin{equation}
\sigma_k\left(H_\lambda\right)=\frac{1}{\lambda}\widetilde{Q}_\lambda(\phi).\nonumber
\end{equation}
It follows that
\begin{equation}
\sigma_k\left(H_\lambda\right)=\int_1^\infty S^{N-1}\left(\upsilon'^2+\frac{1-pu_\lambda^{p-1}}{\lambda}\upsilon^2\right)\,\text{d}r+\mu_{i_k}\int_1^\infty S^{N-3}\upsilon^2\,\text{d}r-(N-1)\frac{e^2+1}{e^2-1},\nonumber
\end{equation}
which is the functional of the following problem
\begin{equation}
\left\{
\begin{array}{lll}
-\lambda\left(\partial_r^2+\frac{(N-1)S'}{S}\partial_r\right)\upsilon+\lambda\frac{\mu_{i_k}}{S^2(r)}\upsilon+\left(1-pu_\lambda^{p-1}\right)\upsilon=0,\\
\upsilon(1)=1.
\end{array}
\right.\nonumber
\end{equation}
By the uniqueness we deduce that $\upsilon(r)\equiv c_{i_k}(r)$.
So we get that
\begin{equation}
\sigma_k\left(H_\lambda\right)=-c_{i_k}'(1)-(N-1)\frac{e^2+1}{e^2-1}=\sigma_{i_k}(\lambda).\nonumber
\end{equation}
Since $E\varsupsetneq E_1^c\varsupsetneq \cdots\varsupsetneq E_k^c\varsupsetneq \cdots$ and $\sigma_k\left(H_\lambda\right)$
is is attained in $E_k$, we have that
\begin{equation}
\sigma_1\left(H_\lambda\right)<\sigma_2\left(H_\lambda\right)<\cdots<\sigma_k\left(H_\lambda\right)<\cdots.\nonumber
\end{equation}
It follows that
\begin{equation}
\sigma_{i_1}\left(\lambda\right)<\sigma_{i_2}\left(\lambda\right)<\cdots<\sigma_{i_k}\left(\lambda\right)<\cdots,\nonumber
\end{equation}
which is the desired conclusion.\qed
\\

We write $\sigma(\lambda)$ as $\sigma_{i_1}(\lambda)$ for simplicity.
Proposition \ref{prop24} implies that $\sigma(\lambda)$ is positive for $\lambda$ large enough
and Proposition \ref{prop25} implies that $\sigma\left(\Lambda_0\right)<0$.
It follows from (\ref{uderivative}) that $u_\lambda$ is continuous with respect to $\lambda$. Then the definitions of $\widetilde{Q}_\lambda$ and ${Q}_\lambda$
imply that $\sigma(\lambda)$ is continuous.
So, $\sigma(\lambda)$ has at least one zero.
We use $\Lambda^*$ to denote the biggest zero of $\sigma(\lambda)$ such that $\sigma$ is negative in a small neighborhood on the left.
Clearly, one see that $\Lambda^*>\Lambda_0$. Take $\Lambda_2>\Lambda^*$ such that $\sigma(\lambda)>0$ for any $\lambda\geq\Lambda_2$.

In view of Proposition \ref{prop43}, we have that
$\sigma_{i_k}\left(\Lambda^*\right)>0$ for any $k>1$. Further we obtain the following result.
\begin{proposition}\label{prop44}
For any $\lambda>\Lambda_0$ and $k\in \mathbb{N}$, $\sigma_{i_k}(\lambda)$ is continuous and possesses at least one zero $\Lambda_k^*$ with $\Lambda^*:=\Lambda_1^*>\Lambda_2^*>\cdots>\Lambda_k^*>\cdots>\Lambda_0$.
\end{proposition}
\textbf{Proof.} The case $k=1$ has been done. For $k>1$, let $\phi(\theta)$ be a normalized eigenfunction corresponding to $\mu_{i_k}$ and $\psi(r,\theta)=u_\lambda(r)\phi(\theta)$.
Similar to Proposition \ref{prop22}, we get that
\begin{eqnarray}
Q_\lambda(\psi)\leq0\nonumber
\end{eqnarray}
for $\lambda\leq -\tau_0S^2(1)/\mu_{i_k}\leq -\tau_0S^2(1)/\mu_{i_1}$.
Then as that of \cite[Lemma 5.4]{RRS} with obvious changes we have that $\sigma_{i_k}\left(\Lambda_0\right)<0$.
Similar to the case of $k=1$, $\sigma_{i_k}(\lambda)$ is also continuous.
Since $\sigma_{i_k}\left(\Lambda^*\right)>0$ for any $k>1$, $\sigma_{i_k}(\lambda)$ has at least one zero. We use $\Lambda_k^*$ the maximum zero of $\sigma_{i_k}$ such that $\sigma_{i_k}$ is negative in the small neighborhood on the left of $\Lambda_k^*$.
The monotonicity of zeros with respect $k$ can be deduced from Proposition \ref{prop43}.\qed



\begin{remark}\label{remark45}
From Proposition \ref{prop44} we see that $\sigma_{i_k}(\lambda)$ with any $k\geq2$ is nonnegative in a small neighborhood on the left of $\Lambda^*$.
Thus, ${H}_\lambda$ has a unique negative eigenvalue in a small neighborhood on the left of $\Lambda^*$.
\end{remark}

\section{Bifurcation and proof the main result}

\quad\, Then similar to \cite[Lemma 6.1]{RRS}, there exists $\varepsilon>0$ such that
for any $\lambda\in\left(\Lambda^*-\varepsilon,+\infty\right)$, the operator $H_\lambda+Id$ is invertible.
\\ \\
\textbf{Proof of Theorem 1.1.}
Take $\Lambda_1\in \left(\Lambda^*-\varepsilon,\Lambda^*\right)$ and define $G:\left[\Lambda_1,\Lambda_2\right]\times \mathcal{V}\rightarrow \mathcal{W}$ by
\begin{equation}
G(\lambda,v)=\frac{F(\lambda,v)}{u_\lambda'(1)}+v,\nonumber
\end{equation}
where $\mathcal{V}\subset C_{G,m}^{2,\alpha}\left(\mathbb{S}^{N-1}\right)$ and $\mathcal{W}\subset C_{G,m}^{1,\alpha}\left(\mathbb{S}^{N-1}\right)$ are open neighborhoods of the $0$ function.
Since the operator $H_\lambda+Id$ is invertible for $\lambda\in \left(\Lambda^*-\varepsilon,+\infty\right)$,
$D_v G(\lambda,0)$ is an isomorphism for all $\lambda\in\left[\Lambda_1,\Lambda_2\right]$. {By using the Inverse Function Theorem, we can further restrict $\mathcal{V}$ and $\mathcal{W}$ so that $G(\lambda, \cdot)$ is invertible for all $\lambda \in [\Lambda_1, \Lambda_2]$.
Let $R(\lambda,w)=w-\tilde w$, where $\tilde w$ is such that $G(\lambda, \tilde w) = w$. Then $R:\left[\Lambda_1,\Lambda_2\right]\times\mathcal{W}\rightarrow
\mathcal{W}$ has the form of identity plus a compact operator.}
Clearly $F(\lambda,v)=0$ is equivalent to $R(\lambda,v)=0$ for all $\lambda\in\left[\Lambda_1,\Lambda_2\right]$.
We see that
$D_w R(\lambda,0)w=\mu w$ is equivalent to $H_\lambda(w)=\mu w/(1-\mu)$ and for $\lambda > \Lambda_1$ we have $\mu<1$.
It follows that $D_w R(\Lambda^*,0)$ and $H_{\Lambda^*}$ have the same kernel space.
Since $m_1$ is odd, the dimension of the kernel space of $D_w R\left(\Lambda^*,0\right)$ is odd.
Then, in view of Remark \ref{remark45}, as that of \cite[Theorem 6.2]{RRS}, applying the Krasnosel'skii Local Bifurcation Theorem \cite{Krasnoselskii} or \cite[Theorem 3.2 in Chapter II]{Kielhofer} to $F(\lambda,v)=0$
, we can conclude the desired bifurcation result.\qed

\bibliographystyle{amsplain}
\makeatletter
\def\@biblabel#1{#1.~}
\makeatother


\begin{thebibliography}{10}




\bibitem{Admas} R.A. Admas, Sobolev spaces, Academic Press, New York, 1975.

\bibitem{Banica} V. Banica, The nonlinear Schr\"{o}dinger equation on hyperbolic space, Comm. Partial Differential Equations 32:10-12
(2007), 1643--1677.

\bibitem{Beckner} W. Beckner, On the Grushin operator and hyperbolic symmetry, Proc. Amer. Math. Soc. 129 (2001), 1233-1246.

\bibitem{Berchio} E. Berchio, A. Ferrero, G. Grillo, Stability and qualitative properties of radial solutions of the Lane-Emden-Fowler
equation on Riemannian models, J. Math. Pures Appl. 102 (1) (2014), 1--35.

\bibitem{BCN} H. Berestycki, L.A. Caffarelli and L. Nirenberg, Monotonicity for elliptic equations in unbounded Lipschitz domains,
Comm. Pure Appl. Math. 50 (1997), 1089--1111.

\bibitem{Born0} M. Born, Modified field equations with a finite radius of the electron, Nature 132 (1933), 282.

\bibitem{Born1} M. Born, Cosmic rays and new field, Nature 133 (1934), 63--64.

\bibitem{Born2} M. Born and L. Infeld, Foundations of the new field theory, Nature 132 (1933), 1004.

\bibitem{Brezis} H. Brezis and L. Nirenberg, Positive solutions of nonlinear elliptic equations involving critial Sobolev exponents, Comm. Pure Appl. Math. 36 (1983), 437-477.

\bibitem{Castorina} D. Castorina, I. Fabbri, G. Mancini and K. Sandeep, Hardy-Sobolev extremals, hyperbolic symmetry and scalar curvature equations, J. Differential Equations, 246 (2009), 1187-1206.

\bibitem{Chavel}  I. Chavel, Eigenvalues in Riemannian geometry, Academic Press, 1984.


\bibitem{DMS} G. Dai, F. Morabito and P. Sicbaldi, A smooth $1$-parameter family of Delaunay-type domains for an overdetermined elliptic problem
in $\mathbb{S}^n\times \mathbb{R}$ and $\mathbb{H}^n\times \mathbb{R}$, Potential Anal., (2023), Preprint.

\bibitem{Del} M. Del Pino, F. Pacard and J. Wei, J. Serrin's overdetermined problem and constant mean curvature
surfaces, Duke Math. J. 164 (2015), 2643--2722.

\bibitem{Espinar} J.M. Espinar, A. Farina and L. Mazet, $f$-extremal domains in hyperbolic space, Calc. Var. Partial Differential Equations 60 (2021), no. 3, Paper No. 112, 19 pp.

\bibitem{EM18} J.M. Espinar and J. Mao, Extremal domains on Hadamard manifolds, J. Differential Equations 265 (2018), 2671--2707.


\bibitem{Fall} M.M. Fall, I.A. Minlend and T. Weth, Unbounded periodic solutions to Serrin's overdetermined boundary value problem, Arch. Ration. Mech. Anal. 223 (2017), 737--759.

\bibitem{Farina} A. Farina and E. Valdinoci, Flattening results for elliptic PDEs in unbounded domains with applications to overdetermined problems, Arch. Ration. Mech. Anal. 195 (2010), 1025--1058.


\bibitem{Ganguly} D. Ganguly and K. Sandeep, Nondegeneracy of positive solutions of semilinear elliptic problems in the hyperbolic space, Commun. Contemp. Math. 17 (2015), no. 1, 1450019, 13 pp.

\bibitem{Gilbarg} D. Gilbarg and N. S. Trudinger, Elliptic partial differential equations of second order, Springer-Verlag, Berlin, Heidelberg, 2001.


\bibitem{Groemer} H. Groemer, Geometric applications of Fourier series and spherical harmonics, Cambridge University Press, 1996.

\bibitem{Helgason} S. Helgason, Geometric analysis on symmetric spaces, Publication Providence, RI:
American Mathematical Society, 1994.

\bibitem{Hofer} H. Hofer, A note on the topological degree at a critical point of mountainpass-type, Proc. Amer. Math. Soc. 90 (1984), 309--315.


\bibitem{Kazdan} J.L. Kazdan, Applications of partial differential equations to problems in geometry, Grad. Texts in Math., Springer, 2004.

\bibitem{Kielhofer} H. Kielhofer, Bifurcation theory: an introduction with applications to PDE's. Appl. Math.
Sci. 156, Springer, New York, 2004.

\bibitem{Krasnoselskii} M.A. Krasnosel'skii, Topological methods in the theory of nonlinear integral equations, Macmillan, New York, 1965.


\bibitem{Kumaresan} S. Kumaresan and J. Prajapat, Serrin's result for hyperbolic space and sphere, Duke Math. J. 91 (1998), 17--28.


\bibitem{Mancini} G. Mancini and K. Sandeep, On a semilinear equation in $\mathbb{H}^n$, Ann. Sc. Norm. Super. Pisa 7 (5) (2008), 635--671.

\bibitem{Morabito0} F. Morabito, Radial and non-radial solutions to an elliptic problem on annular domains in Riemannian manifolds with radial symmetry, J. Differential Equations 258 (2015), 1461--1493.

\bibitem{Morabito} F. Morabito and P. Sicbaldi, Delaunay type domains for an overdetermined elliptic problem in $\mathbb{S}^n\times \mathbb{R}$ and $\mathbb{H}^n\times \mathbb{R}$, ESAIM Control Optim. Calc. Var. 22 (2016), 1--28.


\bibitem{Peletier} L.A. Peletier and J. Serrin, Uniqueness of positive solutions of semilinear equations in $\mathbb{R}^n$, Arch. Rational
Mech. Anal. 81 (1983), 181--197.

\bibitem{Pucci} P. Pucci and J. Serrin, The maximum principle, Progr. Nonlinear differential equations appl., vol. 73, Birkh\"{a}user
Verlag, Basel, 2007.


\bibitem{Reichel} W. Reichel, Radial symmetry for elliptic boundary-value problems on exterior domains, Arch. Ration. Mech.
Anal. 137 (1997), 253--281.

\bibitem{RRS1} A. Ros, D. Ruiz and P. Sicbaldi, A rigidity result for overdetermined elliptic problems in the
plane, Comm. Pure Appl. Math. 70 (2017), 1223--1252.

\bibitem{RRS} A. Ros, D. Ruiz and P. Sicbaldi, Solutions to overdetermined elliptic problems in nontriavial exterior domains, J. Eur. Math. Soc. 22 (2020), 1223--1252.

\bibitem{Ros} A. Ros and P. Sicbaldi, Geometry and topology of some overdetermined elliptic problem, J. Differential Equations 255 (2013), 951--977.

\bibitem{Ruiz} D. Ruiz, Nonsymmetric sign-changing solutions to overdetermined elliptic problems in bounded domains, J. Eur. Math. Soc., to appear.

\bibitem{RSW22} D. Ruiz, P. Sicbaldi, J. Wu. Overdetermined elliptic problems in onduloid-type domains with general nonlinearities, J. Funct. Anal. 283 (2022), 26 pp.

\bibitem{RSW2} D. Ruiz, P. Sicbaldi, J. Wu. Overdetermined elliptic problems in nontrivial contractible domains of the sphere, J. Math. Pures et appl., 180 (2023) 151--187.

\bibitem{Schlenk} F. Schlenk and P. Sicbaldi, Bifurcating extremal domains for the first eigenvalue of the Laplacian, Adv. Math. 229 (2012), 602--632.

\bibitem{Selvitella} A. M. Selvitella, Uniqueness of the ground state of the NLS on $\mathbb{H}^d$ via analytical and topological methods, Rocky Mountain J. Math. 50 (2020), 1817--1832.

\bibitem{Serrin} J. Serrin, A symmetry problem in potential theory, Arch. Ration. Mech. Anal. 43 (1971), 304--318.

\bibitem{Sicbaldi} P. Sicbaldi, New extremal domains for the first eigenvalue of the Laplacian in flat tori, Calc. Var. Partial Differential
Equations 37 (2010), 329--344.




\bibitem{Terras} A. A. Terras, Harmonic analysis on symmetric spaces and applications,
NY-Berlin-Heidelberg: Springer, New York, 1985.

\bibitem{Wang} Z. Wang, Uniqueness of radial solutions of semilinear elliptic equations on hyperbolic space, Nonlinear Anal. 104
(2014), 109--119.


\end{thebibliography}

\providecommand{\bysame}{\leavevmode\hbox to3em{\hrulefill}\thinspace}
\providecommand{\MR}{\relax\ifhmode\unskip\space\fi MR }
\providecommand{\MRhref}[2]{%
  \href{http://www.ams.org/mathscinet-getitem?mr=#1}{#2}
}
\providecommand{\href}[2]{#2}

\end{document}